\documentclass[10pt,a4paper,english]{article}
\usepackage[T1]{fontenc}
\usepackage[latin1]{inputenc}
\usepackage{amsmath}
\usepackage{amssymb}
\usepackage{amsfonts}
\usepackage[all]{xy}
\usepackage{babel}

\providecommand{\tabularnewline}{\\}

\newtheorem{theorem}{Theorem}[section]
\newtheorem{acknowledgement}[theorem]{Acknowledgement}
\newtheorem{algorithm}[theorem]{Algorithm}
\newtheorem{axiom}[theorem]{Axiom}
\newtheorem{case}[theorem]{Case}
\newtheorem{claim}[theorem]{Claim}
\newtheorem{conclusion}[theorem]{Conclusion}
\newtheorem{condition}[theorem]{Condition}
\newtheorem{conjecture}[theorem]{Conjecture}
\newtheorem{corollary}[theorem]{Corollary}
\newtheorem{criterion}[theorem]{Criterion}
\newtheorem{definition}[theorem]{Definition}
\newtheorem{example}[theorem]{Example}
\newtheorem{exercise}[theorem]{Exercise}
\newtheorem{lemma}[theorem]{Lemma}
\newtheorem{notation}[theorem]{Notation}
\newtheorem{problem}[theorem]{Problem}
\newtheorem{proposition}[theorem]{Proposition}
\newtheorem{remark}[theorem]{Remark}
\newtheorem{solution}[theorem]{Solution}
\newtheorem{summary}[theorem]{Summary}
\newenvironment{proof}[1][Proof]{\textbf{#1.} }{\ \rule{0.5em}{0.5em}}

\DeclareMathOperator{\A}{A}
\DeclareMathOperator{\Alb}{Alb}

\DeclareMathOperator{\CH}{CH}

\DeclareMathOperator{\Cp}{Cp}

\DeclareMathOperator{\Decf}{\underline{Dec}}

\DeclareMathOperator{\Div}{Div}

\DeclareMathOperator{\Divf}{\underline{Div}}

\DeclareMathOperator{\LDiv}{LDiv}

\DeclareMathOperator{\WDiv}{WDiv}

\DeclareMathOperator{\Ext}{Ext}
\DeclareMathOperator{\Extf}{\underline{Ext}}

\DeclareMathOperator{\Extfabk}{\underline{Ext}_{\Abk}}

\DeclareMathOperator{\Gl}{GL}

\DeclareMathOperator{\Hom}{Hom}
\DeclareMathOperator{\Homf}{\underline{Hom}}
\DeclareMathOperator{\Homk}{Hom_{\fld}}

\DeclareMathOperator{\Homka}{Hom_{\fld \mathrm{-Alg}}}

\DeclareMathOperator{\Homcka}{Hom_{\fld \mathrm{-Alg}}^{\mathrm{cont}}}
\DeclareMathOperator{\Homabk}{Hom_{\Abk}}

\DeclareMathOperator{\Homfabk}{\underline{Hom}_{\Abk}}

\DeclareMathOperator{\Inf}{Inf}
\DeclareMathOperator{\K}{K}

\DeclareMathOperator{\Res}{Res}

\DeclareMathOperator{\Lie}{Lie}

\DeclareMathOperator{\Mor}{Mor}
\DeclareMathOperator{\Nil}{Nil}

\DeclareMathOperator{\Q}{Q}
\DeclareMathOperator{\R}{R}

\DeclareMathOperator{\Pic}{Pic}

\DeclareMathOperator{\Picf}{\underline{Pic}}

\DeclareMathOperator{\Spec}{Spec}

\DeclareMathOperator{\Spf}{Spf}

\DeclareMathOperator{\Seq}{Seq}
\DeclareMathOperator{\Supp}{Supp}
\DeclareMathOperator{\Sy}{S}
\DeclareMathOperator{\Sym}{Sym}

\DeclareMathOperator{\Z}{Z}

\DeclareMathOperator{\alb}{alb}

\DeclareMathOperator{\cl}{cl}
\DeclareMathOperator{\chr}{char}

\DeclareMathOperator{\dv}{div}
\DeclareMathOperator{\fml}{fml}
\DeclareMathOperator{\hgt}{ht}
\DeclareMathOperator{\id}{id}
\DeclareMathOperator{\im}{im}

\DeclareMathOperator{\inn}{\in}

\DeclareMathOperator{\nil}{nil}

\DeclareMathOperator{\val}{v}
\DeclareMathOperator{\weil}{weil}
\newcommand{\Ab}{\mathbf{Ab}}

\newcommand{\Abk}{\mathcal{A} \mathit{b} / \fld}

\newcommand{\Affk}{\mathbf{Aff} / \fld}

\newcommand{\Algk}{\mathbf{Alg} / \fld}

\newcommand{\Gak}{\mathcal{G} \hspace{-.1em} \mathit{a} / \fld}

\newcommand{\Gfk}{\mathcal{G} \mathit{f} / \fld}

\newcommand{\FSchk}{\mathbf{FSch} / \fld}
\newcommand{\Schk}{\mathbf{Sch} / \fld}
\newcommand{\Set}{\mathbf{Set}}

\newcommand{\Setk}{\mathcal{S} \mathit{et} / \fld}

\newcommand{\Mav}{\mathbf{Mav}}
\newcommand{\Mr}{\mathbf{Mr}}

\newcommand{\MrV}{\mathbf{Vr}}
\newcommand{\MrY}{\mathbf{Yr}}
\newcommand{\Mro}[1]{\Mr_0^{#1}}
\newcommand{\MrCH}[1]{\Mr^{\CHOo {#1}}}

\newcommand{\Divoyx}{\Divf^{0}_{Y/X}}

\newcommand{\Homabs}[1]{\mathrm{Hom}_{\mathcal{A} \mathit{b} / #1}}

\renewcommand{\H}{\mathrm{H}}

\newcommand{\Insg}{\mathrm{H}_{\mathrm{I}}^{-1} \,}

\newcommand{\CHO}[1]{\CH_0 \left(#1\right)}
\newcommand{\CHOo}[1]{\CH_0 \left(#1\right)^{0}}

\newcommand{\ZOo}[1]{\Z_0 \left(#1\right)^{0}}

\newcommand{\Nat}{\mathbb{N}}
\newcommand{\Zint}{\mathbb{Z}}

\newcommand{\Cplx}{\mathbb{C}}

\newcommand{\Prj}{\mathbb{P}}

\newcommand{\cont}{\mathrm{cont}}
\newcommand{\cst}{\mathrm{const}}
\newcommand{\der}{\mathrm{d}}

\newcommand{\et}{\acute{\mathrm{e}} \mathrm{t}}

\newcommand{\htl}{\hgt = \mathrm{1}}

\newcommand{\lib}{\mathrm{lib}}
\newcommand{\mdl}{\mathfrak{d}}
\newcommand{\mdll}{\mathfrak{e}}

\newcommand{\ord}{\mathrm{ord}}

\newcommand{\reg}{\mathrm{reg}}

\newcommand{\sing}{\mathrm{sing}}
\newcommand{\tor}{\mathrm{tor}}
\newcommand{\trafo}{\tau}
\newcommand{\transl}{\mathrm{t}}

\newcommand{\ke}{k[\varepsilon]}

\newcommand{\fld}{k}

\newcommand{\fmlE}{\mathcal{E}}
\newcommand{\fmlF}{\mathcal{F}}
\newcommand{\fmlG}{\mathcal{F}}
\newcommand{\fmlGr}{\mathcal{G}}

\newcommand{\splG}{\mathbb{G}}

\newcommand{\Trs}{T}
\newcommand{\Trsd}{\Trs^{\vee}}

\newcommand{\Vcl}{\mathbb{V}}
\newcommand{\Vcld}{\Vcl^{\vee}}

\newcommand{\Ga}{\mathbb{G}_{\mathrm{a}}}

\newcommand{\GaS}[1]{\mathbb{G}_{\mathrm{a},#1}}
\newcommand{\Gm}{\mathbb{G}_{\mathrm{m}}}
\newcommand{\GmS}[1]{\mathbb{G}_{\mathrm{m},#1}}
\newcommand{\pnt}{q}
\newcommand{\pntt}{p}

\newcommand{\Es}{E}

\newcommand{\Sing}{S}
\newcommand{\SgY}{\Sing_Y}
\newcommand{\Curv}{Z}
\newcommand{\Crvv}{C}

\newcommand{\pt}{\widetilde{p}}

\newcommand{\mc}{\widehat{\fm}}

\newcommand{\Oc}{\widehat{\sO}}

\newcommand{\Vt}{\widetilde{V}}

\newcommand{\Xt}{\widetilde{X}}
\newcommand{\Ct}{\widetilde{C}}

\newcommand{\Ld}{L^{\vee}}
\newcommand{\Ad}{A^{\vee}}
\newcommand{\Gd}{G^{\vee}}
\newcommand{\Fd}{\sF^{\vee}}
\newcommand{\ld}{l^{\vee}}

\newcommand{\sigt}{\widetilde{\sigma}}

\newcommand{\wh}[1]{\widehat{#1}}

\newcommand{\lp}{\left(}
\newcommand{\rp}{\right)}
\newcommand{\ra}{\rightarrow}

\newcommand{\dra}{\dashrightarrow}

\newcommand{\lra}{\longrightarrow}
\newcommand{\lla}{\longleftarrow}

\newcommand{\Lra}{\Longrightarrow}
\newcommand{\Lla}{\Longleftarrow}
\newcommand{\Llra}{\Longleftrightarrow}
\newcommand{\lmt}{\longmapsto}
\newcommand{\mpt}{\mapsto}

\newcommand{\inj}{\hookrightarrow}
\newcommand{\sur}{\twoheadrightarrow}
 
\newcommand{\cut}{\cdot}
\newcommand{\isec}{\cap}
\newcommand{\dsum}{\bigoplus}

\newcommand{\tens}{\otimes}

\newcommand{\tms}{\times}

\newcommand{\fis}{\#}

\newcommand{\see}{}

\newcommand{\alp}{\alpha}
\newcommand{\gam}{\gamma}
\newcommand{\del}{\delta}
\newcommand{\eps}{\varepsilon}
\newcommand{\tha}{\vartheta}

\newcommand{\sig}{\sigma}

\newcommand{\phe}{\varphi}
\newcommand{\oma}{\omega}

\newcommand{\Gam}{\Gamma}

\newcommand{\Oma}{\Omega}

\newcommand{\fm}{\mathfrak{m}}

\newcommand{\fR}{\mathfrak{R}}

\newcommand{\sA}{\mathcal{A}}

\newcommand{\sD}{\mathcal{D}}

\newcommand{\sF}{\mathcal{F}}
\newcommand{\sG}{\mathcal{G}}

\newcommand{\sK}{\mathcal{K}}
\newcommand{\sL}{\mathcal{L}}

\newcommand{\sO}{\mathcal{O}}

\newcommand{\sU}{\mathcal{U}}
\newcommand{\sV}{\mathcal{V}}

\newcommand{\sY}{\mathcal{Y}}

\newcommand{\bThm}{\begin{theorem}}
\newcommand{\eThm}{\end{theorem}}
\newcommand{\bAck}{\begin{acknowledgement}}
\newcommand{\eAck}{\end{acknowledgement}}
\newcommand{\bAlg}{\begin{algorithm}}
\newcommand{\eAlg}{\end{algorithm}}
\newcommand{\bAxm}{\begin{axiom}}
\newcommand{\eAxm}{\end{axiom}}
\newcommand{\bCas}{\begin{case}}
\newcommand{\eCas}{\end{case}}
\newcommand{\bClm}{\begin{claim}}
\newcommand{\eClm}{\end{claim}}
\newcommand{\bCcl}{\begin{conclusion}}
\newcommand{\eCcl}{\end{conclusion}}
\newcommand{\bCdn}{\begin{condition}}
\newcommand{\eCdn}{\end{condition}}
\newcommand{\bCjc}{\begin{conjecture}}
\newcommand{\eCjc}{\end{conjecture}}
\newcommand{\bCor}{\begin{corollary}}
\newcommand{\eCor}{\end{corollary}}
\newcommand{\bCrt}{\begin{criterion}}
\newcommand{\eCrt}{\end{criterion}}
\newcommand{\bDef}{\begin{definition}}
\newcommand{\eDef}{\end{definition}}
\newcommand{\bExm}{\begin{example}}
\newcommand{\eExm}{\end{example}}
\newcommand{\bExc}{\begin{exercise}}
\newcommand{\eExc}{\end{exercise}}
\newcommand{\bLem}{\begin{lemma}}
\newcommand{\eLem}{\end{lemma}}
\newcommand{\bNot}{\begin{notation}}
\newcommand{\eNot}{\end{notation}}
\newcommand{\bPrb}{\begin{problem}}
\newcommand{\ePrb}{\end{problem}}
\newcommand{\bPrp}{\begin{proposition}}
\newcommand{\ePrp}{\end{proposition}}
\newcommand{\bRmk}{\begin{remark}}
\newcommand{\eRmk}{\end{remark}}
\newcommand{\bSol}{\begin{solution}}
\newcommand{\eSol}{\end{solution}}
\newcommand{\bSmr}{\begin{summary}}
\newcommand{\eSmr}{\end{summary}}
\newcommand{\bPf }{\begin{proof}}
\newcommand{\ePf }{\end{proof}}

\newcommand{\vs}{12pt}

\begin{document}

\centerline{ } 
\vspace{25pt} 
\centerline{\LARGE{Generalized Albanese and its dual}}
\vspace{20pt}
\centerline{\normalsize{Henrik Russell}\footnote{
This work was supported by the DFG Leibniz prize 
and the DFG Graduiertenkolleg 647}
} 
\vspace{10pt}
\centerline{\normalsize{August 2008}} 
\vspace{20pt}

\begin{abstract} 
Let $X$ be a projective variety over an algebraically closed field $k$ 
of characteristic 0. 
We consider categories of rational maps from $X$ to commutative 
algebraic groups, and ask for objects satisfying the universal 
mapping property. 
A necessary and sufficient condition for the existence of such 
universal objects is given, as well as their explicit construction, 
using duality theory of generalized 1-motives. 

An important application is the Albanese of a singular projective 
variety, which was constructed by Esnault, Srinivas and Viehweg 
as a universal regular quotient of a relative Chow group of 
0-cycles of degree 0 modulo rational equivalence. 
We obtain functorial descriptions of the universal regular quotient 
and its dual 1-motive. 
\end{abstract}

\tableofcontents 
\setcounter{section}{-1}

%\newpage 

\section{Introduction} 

For a projective variety $X$ over an algebraically closed field $k$ 
a generalized Albanese variety $\Alb\left(X\right)$ is constructed
by Esnault, Srinivas and Viehweg in \cite{ESV} 
as a universal regular quotient of the relative Chow-group $\CHOo{X}$ 
of Levine-Weibel \cite{LW} 
of $0$-cycles of degree $0$ modulo rational equivalence. 
This is a smooth connected commutative algebraic group, 
universal for rational maps from $X$ to smooth commutative 
algebraic groups $G$ factoring through a homomorphism of groups 
$\CHOo{X}\lra G(k)$. 
It is not in general an abelian variety if $X$ is singular. 
Therefore it cannot be dualized in the same way as an abelian variety. 

Laumon built up in \cite{L} a duality theory of generalized 1-motives
in characteristic $0$, which are homomorphisms $\left[\sF\ra G\right]$
from a commutative torsion-free formal group $\sF$ to a connected 
commutative algebraic group $G$. 
The universal regular quotient $\Alb\left(X\right)$ can be interpreted 
as a generalized 1-motive by setting $\sF=0$ and $G=\Alb\left(X\right)$. 
The motivation for this work was to find the functor 
which is represented by the dual 1-motive, 
in the situation where the base field $k$ is 
algebraically closed and of characteristic 0. 
The duality gives an independent proof 
(alternative to the ones in \cite{ESV}) 
as well as an explicit construction of the universal regular quotient 
in this situation 
(cf.\ Subsection \ref{subsec:Universal_Regular_Quotient}). 
This %(Theorem \ref{main_result}) 
forms one of the two main results of the present article: 

\bThm \label{main_result} 
Let $X$ be a projective variety over an algebraically closed field $k$ 
of characteristic 0, and $\Xt\lra X$ a projective resolution of singularities. 
Then the universal regular quotient $\Alb\left(X\right)$ exists and 
there is a subfunctor $\Divf_{\Xt/X}^0$ of the relative Cartier divisors 
on $\Xt$ (cf.\ Subsection \ref{sub: Div_Y/X^0}, Definition \ref{Div_Y/X^0}) 
such that the dual of $\Alb\left(X\right)$ 
(in the sense of 1-motives) represents the functor 
\[ \Divf_{\Xt/X}^{0} \lra \Picf_{\Xt}^{0}
\] 
i.e.\ the natural transformation of functors 
which assigns to a relative Cartier divisor 
the class of its associated line bundle. 
$\Picf_{\Xt}^{0}$ is represented by an abelian variety and
$\Divf_{\Xt/X}^{0}$ by a formal group. 
\eThm 

The other main result (Theorem \ref{univ_object}) is a more general statement 
about the existence and construction of universal objects 
of categories of rational maps 
(the notion of \emph{category of rational maps} 
is introduced in Definition \ref{CatMr}, but the name is suggestive). 
This concept does not only contain the universal regular quotient, 
but also the generalized Jacobian of Rosenlicht-Serre \cite[Chapter V]{S} 
as well as the generalized Albanese of Faltings-W\"ustholz \cite{FW} 
as special cases of such universal objects. 

%One might ask why 
We only deal with a base field of characteristic 0, 
although the universal regular quotient exists in any characteristic. 
A first reason for this is 
that Laumon's 1-motives are only defined in characteristic 0. 
In order to match the case of arbitrary characteristic, one first needs 
to define a new category of 1-motives (in any characteristic) 
which contains smooth connected commutative algebraic groups 
as a subcategory. 
A commutative torsion-free 
formal group in characteristic 0 
is completely determined by its $k$-valued points and its Lie-algebra 
(cf.\ Corollary \ref{fmlG determin}), 
the first form a free 
abelian group of finite rank, 
the latter is a finite dimensional $k$-vector space. 
%This is possible by means of the exponential map. 
This allows to give an explicit and transparent description, 
which might not be possible in positive characteristic.

\subsection{Leitfaden} 

In the following we give a short summary of each section. 

\textbf{Section \ref{sec:1-Motives}} provides some basic facts about
generalized 1-motives, which are used in the rest of the paper. 
A connected commutative algebraic group $G$ is an extension 
of an abelian variety $A$ by a linear group $L$. 
Then the dual 1-motive of $\left[0\ra G\right]$
is given by $\left[\Ld\ra\Ad\right]$, 
where $\Ld=\Homf\left(L,\Gm\right)$ is the Cartier-dual of $L$ and 
$\Ad=\Pic_{A}^{0}=\Extf\left(A,\Gm\right)$ is the dual abelian variety, 
and the homomorphism between them is the connecting homomorphism in the 
long exact cohomology sequence obtained from $0 \ra L \ra G \ra A \ra 0$ 
by applying the functor $\Homf\left(\_,\Gm\right)$. 

\textbf{Section \ref{sec:Univ-Fact-Prbl}} states the universal factorization
problem with respect to a category $\Mr$ of rational maps from a
regular projective variety $Y$ to connected commutative algebraic
groups (cf.\ Definition \ref{UnivObj}): 
\bDef  
A rational map $\left(u:Y\dra\sU\right)\in\Mr$
is called \emph{universal} for $\Mr$ if for all objects 
$\left(\phe:Y\dra G\right)\in\Mr$
there is a unique homomorphism of algebraic groups $h:\sU\lra G$
such that $\phe=h\circ u$ up to translation. 
\eDef 

An essential ingredient for the construction of such universal objects is 
the functor of relative Cartier divisors $\Divf_{Y}$ on $Y$, 
which assigns to an affine scheme $T$ a family of Cartier divisors on $Y$, 
parameterized by $T$. 
This functor admits a natural transformation $\cl:\Divf_Y\lra\Picf_Y$ to 
the Picard functor $\Picf_{Y}$, which maps a relative divisor to its class. 
Then $ \Divf_{Y}^{0}:=\cl^{-1}\Picf_{Y}^{0} $ 
is the functor of families of Cartier divisors whose associated line
bundles are algebraically equivalent to the trivial bundle. 

We give a necessary and sufficient condition for the existence of
a universal object for a category of rational maps $\Mr$ which contains
the category $\Mav$ of morphisms from $Y$ to abelian varieties
and satisfies a certain stability condition $\left(\diamondsuit\right)$, 
see Subsection \ref{sub:Universal-Objects}. 
%$\see$page~\pageref{Exist univObj}. 
Localization of $\Mr$ at the system of injective homomorphisms 
does not change the universal object; 
denote this localization by $\Insg\Mr$. 
We observe that a rational map $\phe:Y\dra G$, where
$G$ is an extension of an abelian variety by a linear group $L$,
induces a natural transformation $\Ld\lra\Divf_{Y}^{0}$.
If $\fmlG$ is a formal group that is a subfunctor of $\Divf_{Y}^{0}$, 
then $\Mr_{\fmlG}$ denotes the category of rational maps for which the
image of this induced transformation lies in $\fmlG$. 
We show (cf.\ Theorem \ref{Exist univObj}):  

\bThm \label{univ_object}
For a category $\Mr$ containing $\Mav$ and satisfying 
$\left(\diamondsuit\right)$ there exists a universal object 
$\Alb_{\Mr}\left(Y\right)$ if and only if there is a formal group 
$\fmlG\subset\Divf_{Y}^{0}$ such that $\Insg\Mr$ is equivalent to 
$\Insg\Mr_{\fmlG}$. 
\eThm

The universal object $\Alb_{\fmlG}\left(Y\right)$ of $\Mr_{\fmlG}$ is an 
extension of the classical Albanese $\Alb\left(Y\right)$,
which is the universal object of $\Mav$, by the linear group $\fmlG^{\vee}$,
the Cartier-dual of $\fmlG$. The dual 1-motive of 
$\left[0\lra\Alb_{\fmlG}\left(Y\right)\right]$
is hence given by $\left[\fmlG\lra\Pic_{Y}^{0}\right]$, 
which is the homomorphism induced by the natural transformation 
$\cl:\Divf_{Y}^{0}\lra\Picf_{Y}^{0}$. 

The universal regular quotient $\Alb\left(X\right)$ of a (singular)
projective variety $X$ is by definition the universal object for
the category $\MrCH{X}$ of rational maps factoring through rational 
equivalence. 
More precisely, the objects of $\MrCH{X}$ are rational maps $\phe:X\dra G$ 
whose associated map on zero-cycles of degree zero 
$ 
  \ZOo{U} \lra G(k), \;
  \sum n_i \, p_i \lmt \sum n_i \, \phe(p_i) 
$ 
(here $U$ is the open set on which $\phe$ is defined) 
factors through a homomorphism of groups $\CHOo{X} \lra G(k)$, 
where $\CHO{X}$ denotes the relative Chow group of 0-cycles 
$\CH_0(X,X_{\sing})$ in the sense of \cite{LW}. 
Such a rational map is regular on the regular locus of $X$ and may also be 
considered as a rational map from $\Xt$ to $G$, where $\pi:\Xt\lra X$ 
is a projective resolution of singularities. 
In particular, if $X$ is nonsingular, 
the universal regular quotient coincides with the classical Albanese. 
The category $\MrCH{X}$ contains $\Mav$ and satisfies 
$\left(\diamondsuit\right)$. 
Therefore our problem reduces to finding the subfunctor of $\Divf_{\Xt}^{0}$ 
which is represented by a formal group $\fmlG$ such that 
$\Mr_{\fmlG}$ equivalent to $\MrCH{X}$. 

\textbf{Section \ref{sec:Rat-Maps CH_0(X)_deg0}} answers the question
for the formal group $\fmlG$ which characterizes the category 
$\MrCH{X}$.
This is a subfunctor of $\Divf_{\Xt}^{0}$ which measures the difference 
between $\Xt$ and $X$. 

If $\pi:Y\lra X$ is a proper birational morphism of varieties, 
the \emph{push-forward} of cycles gives a homomorphism 
$\pi_{*}:\WDiv(Y)\lra\WDiv(X)$ 
from the group of Weil divisors on $Y$ to the group of those on $X$.  
For a curve $\Crvv$ we introduce the $k$-vector space $\LDiv(\Crvv)$ 
of formal infinitesimal divisors,
which generalizes infinitesimal deformations of the zero divisor. 
There exist non-trivial infinitesimal deformations of Cartier divisors, 
but not of Weil divisors, since prime Weil divisors are always reduced. 
Formal infinitesimal divisors also admit a 
\emph{push-forward} $\pi_{*}:\LDiv(\Curv)\lra\LDiv(\Crvv)$ 
for finite morphisms $\pi:\Curv\lra\Crvv$ of curves, 
e.g.\ for the normalization, which is a resolution of singularities. 
There exist natural homomorphisms \;$\weil:\Divf_Y(k)\lra\WDiv(Y)$ 
and \;$\fml:\Lie(\Divf_{\Curv})\lra\LDiv(\Curv)$. 

For a curve $\Crvv$, a natural candidate for the formal group 
we are looking for is determined by the following conditions: 
\begin{eqnarray*}
\Divf_{\Ct/\Crvv}^{0}(k) & = & 
\ker\left(\Divf_{\Ct}^{0}(k)\overset{\weil}\lra
    \WDiv\big(\Ct\big)\overset{\pi_{*}}\lra\WDiv(\Crvv)\right) \\ 
\Lie\left(\Divf_{\Ct/\Crvv}^{0}\right) & = & 
\ker\left(\Lie\left(\Divf_{\Ct}^{0}\right)\overset{\fml}\lra
    \LDiv\big(\Ct\big)\overset{\pi_{*}}\lra\LDiv(\Crvv)\right). 
\end{eqnarray*} 
For a higher dimensional variety $X$, the definition is derived 
from the one for curves as follows. 
A morphism of varieties $V\lra Y$ induces a natural transformation: 
the \emph{pull-back} of relative Cartier divisors 
$\_\cut V:\Decf_{Y,V}\lra\Divf_{V}$, 
where $\Decf_{Y,V}$ is the subfunctor of $\Divf_Y$ consisting of 
those relative Cartier divisors on $Y$ which do not contain 
any component of $\im(V\lra Y)$. 
Then we let $\Divf_{\Xt/X}^0$ be the formal subgroup of $\Divf_{\Xt}^0$ 
characterized by the conditions 
\begin{eqnarray*}
    \Divf_{\Xt/X}^0(k) & = & 
    \bigcap_{C} \left( \_\cut\Ct \right)^{-1} \Divf_{\Ct/C}^0(k) \\
    \Lie\big(\Divf_{\Xt/X}^0\big) & = & 
    \bigcap_{C} \left( \_\cut\Ct \right)^{-1} \Lie\big(\Divf_{\Ct/C}^0\big) \\
\end{eqnarray*} 
where the intersection ranges over all Cartier curves in $X$ relative to 
the singular locus of $X$ in the sense of \cite{LW}. 
Actually, it is not necessary to consider all Cartier curves, 
the functor $\Divf_{\Xt/X}^0$ can be computed from one single general curve. 
%of a suitable smaller family of Cartier curves in $X$. 
The verification of the equivalence between $\MrCH{X}$
and $\Mr_{\Divf_{\Xt/X}^{0}}$ is done using \emph{local symbols}, 
for which \cite{S} is a good reference. 

This gives an independent proof of the existence 
(alternative to the ones in \cite{ESV}) 
as well as an explicit construction of the universal regular quotient 
over an algebraically closed base field of characteristic 0 
(cf.\ Subsection \ref{subsec:Universal_Regular_Quotient}). 

The universal regular quotient for semi-abelian varieties, i.e.\ the universal 
object for rational maps to semi-abelian varieties factoring through rational 
equivalence (which is a quotient of our universal regular quotient), 
is a classical 1-motive in the sense of Deligne 
\cite[D\'efinition (10.1.2)]{D2}. 
The question for the dual 1-motive of this object was already 
answered by Barbieri-Viale and Srinivas in \cite{BS}. 

\textbf{Acknowledgement.} 
This article constitutes the heart of my PhD thesis. 
I am very grateful 
to my advisors H\'el\`ene Esnault and Eckart Viehweg for their support 
and guidance, and for the interesting subject of my thesis, which 
fascinated me from the first moment. 
I owe thanks to Kay R\"ulling for an important hint. 
Moreover, I would like to thank Kazuya Kato for his interest in my studies 
and many helpful discussions.

\subsection{Notations and Conventions}
\label{Terminology} 

$k$ is a fixed algebraically closed field of characteristic $0$. 
%Schemes are always locally Noetherian. 
A variety is a reduced scheme of finite type over $k$. 
%not necessarily irreducible. 
A curve is a variety of dimension 1. 
Algebraic groups and formal groups are always commutative and over $k$. 
We write $\Ga$ for $\GaS{k}=\Spec k[t]$ and $\Gm$ for 
$\GmS{k}=\Spec k\left[t,t^{-1}\right]$. 
The letter $\splG$ stands for a linear algebraic group which is either 
$\Ga$ or $\Gm$. 

If $Y$ is a scheme, then $y\in Y$ means that $y$ is a point in the 
Zariski topological space of $Y$. 
The set of irreducible components of $Y$ is denoted by $\Cp(Y)$. 
%the set of connected components by $\CCp(Y)$. 

If $M$ is a module, then $\Sym M$ is the symmetric algebra of $M$. 
If $A$ is a ring, then $\K_{A}$ denotes the total quotient ring
of $A$. If $Y$ is a scheme, then $\sK_{Y}$ denotes the sheaf
of total quotient rings of $\sO_{Y}$. The group of units of a ring 
$R$ is denoted by $R^*$. 

If $\sig: Y \lra X$ is a morphism of schemes, then 
$\sig^{\fis}: \sO_X \lra \sig_*\sO_Y$ denotes the associated 
homomorphism of structure sheaves. 
%If $h: A \lra B$ is a homomorphism of rings, then 
%$h^{\cis}: \Spec B \lra \Spec A$ denotes the associated 
%morphism of affine schemes. 

We think of $\Ext^{1}\left(A,B\right)$ as the space of extensions
of $A$ by $B$ and therefore denote it by $\Ext\left(A,B\right)$.

The dual of an object $O$ in its respective category is denoted by $O^{\vee}$, 
whereas $\widehat{O}$ is the completion of $O$. 
For example, if $V$ is a $k$-vector space,
then $V^{\vee}=\Homk\left(V,k\right)$ is the dual $k$-vector space; 
if $G$ is a linear algebraic group or a formal group, 
then $G^{\vee}=\Homfabk\left(G,\Gm\right)$ is the Cartier-dual; 
if $A$ is an abelian variety, 
then $A^{\vee}=\Pic^0 A$ is the dual abelian variety. 
%whereas $\widehat{A}=\Spf\widehat{\sO}_{A,0}$ is the 
%completion of $A$ w.r.t. $0_{A}$. 
 
%\newpage 

\section{1-Motives} 
\label{sec:1-Motives} 

The aim of this section is to summarize some foundational material
about generalized 1-motives (following \cite[Sections 4 and 5]{L}), 
as far as necessary for the purpose of these notes. \\ 
Throughout the whole work the base field $k$ is algebraically closed 
and of characteristic 0.

\subsection{Algebraic Groups and Formal Groups }

Here we recall some basic facts about algebraic groups
and the notion of a formal group. 
References for formal groups and Cartier duality 
are e.g.\ \cite[$\mathrm{VII}_\mathrm{B}$]{SGA3}, 
\cite[Chapter II]{D} and \cite[Chapitre I]{Fo}.

\subsubsection*{Algebraic Groups} 
\label{sub:Algebraic-Groups} 

An \emph{algebraic group} %(or \emph{group scheme}) 
is a commutative group-object in the category of 
separated schemes of finite type over $k$. 
As $\chr\left(k\right)=0$, an algebraic group is always smooth 
(see \cite[Chapter III, No.~11, p.~101]{M}). 

\bThm[Chevalley] \label{Thm Chevalley}
A smooth connected algebraic group $G$ admits a \emph{canonical decomposition}
\[ 0\lra L\lra G\lra A\lra0 
\] 
where $L$ is a connected linear algebraic group and $A$ is an abelian variety.
\eThm 
(See \cite[Section 5, Theorem 16, p.~439]{R} or 
\cite[Theorem 3.2, p.~97]{B} or \cite{C}.) 

\bThm  \label{Thm linGrp} A connected linear algebraic group $L$ 
splits canonically into a direct product of a torus $\Trs$ 
and a unipotent group $\Vcl$:
\[ L = \Trs \times_k \Vcl \;. 
\] 
A torus over an algebraically closed field is the direct product 
of several copies of the multiplicative group $\Gm$. 
For $\chr\left(k\right)=0$ a unipotent group is always vectorial, i.e.\ 
is the direct product of several copies of the additive group $\Ga$. 
\eThm 
(See \cite[Expos\'e $\textrm{XVII}$, 7.2.1]{SGA3} and \cite[(4.1)]{L}.)

\subsubsection*{Formal Groups }
\label{sub:Formal-Groups}

A \emph{formal scheme} over $k$ is a functor 
from the category of $k$-algebras $\Algk$ to the category of sets $\Set$ 
which is the inductive limit of a directed inductive set 
of finite $k$-schemes: 
$\sF$ is a formal scheme if there exists a directed projective system $(A_i)$ 
of finite dimensional $k$-algebras and an isomorphism of functors 
$\sF \cong \varinjlim \Spf A_i$, 
where $\Spf A_i: \Algk \lra \Set$ is the functor $R \lmt \Homka(A_i,R)$. 
Equivalently, there is a profinite $k$-algebra $\sA$, 
i.e.\ $\sA$ is the projective limit (as a topological ring) 
of discrete quotients 
which are finite dimensional $k$-algebras, 
and an isomorphism of functors $\sF \cong \Spf\sA$, 
where $\Spf\sA$ is the functor which assigns to a $k$-algebra $R$ 
the set of continuous homomorphisms of $k$-algebras from $\sA$ 
to the discrete ring $R$. 
(Cf. \cite[Chapter I, No.~6]{D}.) 

A \emph{formal group} $\sG$ is a commutative group-object
in the category of formal schemes over $k$, 
such that $\sG(k)$ is an abelian group of finite type 
and $\fm_{\sG,0}/\fm_{\sG,0}^2$ is a finite dimensional $k$-vector space, 
where $\fm_{\sG,0}$ is the maximal ideal of the local ring $\sO_{\sG,0}$. 
If $\chr\left(k\right)=0$, a formal group $\sG$ is always 
equi-dimensional and formally smooth, 
i.e.\ there is a natural number $d\geq0$ such that 
$\sO_{\sG,0}\cong k\left[\left[x_{1},\ldots,x_{d}\right]\right]$ 
(cf.\ \cite[(4.2)]{L}). 

\bThm \label{Str-fml-Grp} 
A formal group $\sG$ admits a canonical decomposition 
\[ \sG \,\cong\, \sG_{\et}\times\sG_{\inf} \] 
where $\sG_{\et}$ is \'etale over $k$ and $\sG_{\inf}$ 
is the component of the identity 
(called \emph{infinitesimal} formal group). 
\eThm 
(See \cite[Chapitre I, 7.2]{Fo} or \cite[(4.2.1)]{L}.) 

\bThm  \label{G_et(R)}
An \'etale formal group $\sG_{\et}$ admits a canonical decomposition 
\[ 0\lra\sG_{\et}^{\tor}\lra\sG_{\et}\lra\sG_{\et}^{\lib}\lra0 
\] 
where $\sG_{\et}^{\tor}$ is the largest sub-group scheme whose underlying
$k$-scheme is finite and \'etale, and $\sG_{\et}^{\lib}\left(k\right)$
is a free abelian group of finite rank. 
\eThm 
(See \cite[(4.2.1)]{L}.) 

\bThm \label{ifGr-Lie} 
For $\chr\left(k\right)=0$, the Lie-functor gives 
an equivalence between the following categories: 
\[ \left\{ \textrm{infinitesimal formal groups}/k\right\} 
   \longleftrightarrow\left\{ \textrm{finite dim. vector spaces}/k\right\} 
\]
\eThm 
(Cf. \cite[$\textrm{VII}_{\textrm{B}}$, 3.3.2.]{SGA3}.)

%\vspace{\vs} 

\bRmk \label{G_inf(R)} 
Theorem \ref{ifGr-Lie} says that 
for an infinitesimal formal group $\sG_{\inf}$, 
there is a finite dimensional $k$-vector space $V$, 
namely $V=\Lie\left(\sG_{\inf}\right)$, such that 
$\sG_{\inf}\cong\Spf\left(\widehat{\Sym V^{\vee}}\right)$, 
where $\wh{\Sym V^{\vee}}$ is the completion of the symmetric algebra 
$\Sym V^{\vee}$ w.r.t. the ideal generated by $V^{\vee}$. 
Therefore the $R$-valued points of $\sG_{\inf}$ 
are given by $\sG_{\inf}(R) = 
\Homcka\left(
  \widehat{\Sym\left(\Lie\left(\sG_{\inf}\right)^{\vee}\right)},R\right) 
= \Lie\left(\sG_{\inf}\right) \tens_k \Nil(R)$. 
\eRmk 
{}From these structure theorems we obtain 

\bCor \label{fmlG determin} 
A formal group $\sG$ in characteristic 0 is uniquely determined by its 
$k$-valued points $\sG(k)$ and its Lie-algebra $\Lie\left(\sG\right)$. 
\eCor

\subsubsection*{Sheaves of Abelian Groups}

The category of algebraic groups and the category of formal groups
can be viewed as full subcategories of the category of sheaves
of abelian groups: 

\bDef Let 

\begin{tabular}{ll}
$\Ab$&
category of abelian groups, \tabularnewline
$\Algk$&
category of $k$-algebras, \tabularnewline
$\Affk$&
category of affine $k$-schemes, \tabularnewline
$\Schk$&
category of $k$-schemes, \tabularnewline
$\FSchk$&
category of formal $k$-schemes. \tabularnewline
\end{tabular} \\
$\Affk$ is anti-equivalent to $\Algk$. Let $\Affk$ (resp. $\Algk$) 
be equipped with the topology fppf. 
Let 

\begin{tabular}{ll}
$\Setk$&
category of sheaves of sets over $\Affk$, \tabularnewline
$\Abk$&
category of sheaves of abelian groups over $\Affk$, \tabularnewline
$\Gak$&
category of algebraic groups over $k$, \tabularnewline
$\Gfk$&
category of formal groups over $k$. \tabularnewline
\end{tabular}
\eDef 
Interpreting a $k$-scheme $X$ as a sheaf over $\Affk$ given by
\[ S\lmt X\left(S\right)=\Mor_{k}\left(S,X\right) \] 
or equivalently over $\Algk$ 
\[ R\lmt X\left(R\right)=\Mor_{k}\left(\Spec R,X\right) \]
makes $\Schk$ a full subcategory of $\Setk$ and $\Gak$ 
a full subcategory of $\Abk$. 

In the same manner $\FSchk$ becomes a full subcategory of $\Setk$
and $\Gfk$ a full subcategory of $\Abk$: 
A formal $k$-scheme $\sY=\Spf\sA$, where $\sA$ is a profinite $k$-algebra,
is viewed as the sheaf over $\Affk$ given by \[
R\lmt\sY\left(R\right)=\underline{\Spf\sA}\left(R\right)=
\Homcka\left(\sA,R\right)\]
which assigns to a $k$-algebra $R$ with discrete
topology the set of continuous homomorphisms of $k$-algebras from
$\sA$ to $R$. 

The categories $\Gak$ and $\Gfk$ are abelian 
(see \cite[(4.1.1) and (4.2.1)]{L}). 
Kernel and cokernel of a homomorphism in $\Gak$ (resp. $\Gfk$) 
coincide with the ones in $\Abk$, and an exact sequence
$0\ra K\ra G\ra C\ra0$ in $\Abk$, 
where $K$ and $C$ are objects of $\Gak$ (resp. $\Gfk$), 
implies that $G$ is also an object of $\Gak$ (resp. $\Gfk$).

\subsection{Definition of 1-Motive} 

In the following by a 1-motive always a \emph{generalized 1-motive}
in the sense of Laumon \cite[D\'efinition (5.1.1)]{L} is meant: 

\bDef 
\label{1-motive} 
A \emph{1-motive} is a complex concentrated in degrees $-1$
and $0$ in the category of sheaves of abelian groups of the form
$M=\left[\sF\ra G\right]$, where $\sF$ is a torsion-free
formal group over $k$ and $G$ a connected algebraic group over $k$.
\eDef

\subsection{Cartier-Duality} 

Let $G$ be an algebraic or a formal group and let $\Homfabk\left(G,\Gm\right)$
be the sheaf of abelian groups over $\Algk$ associated to the functor
\[ R \lmt \Homabs{R}\left(G_R,\GmS{R}\right) 
\] 
which assigns to a $k$-algebra $R$ the set of 
homomorphisms of sheaves of abelian groups over $R$ 
from $G_R$ to $\GmS{R}$. 
If $G$ is a linear algebraic group (resp. formal group), this functor 
is represented by  a formal group (resp. linear algebraic group) $\Gd$, 
called the \emph{Cartier-dual} of $G$. 

The Cartier-duality is an anti-equivalence between the category
of linear algebraic groups and the category of formal groups.
The functors $L\lmt\Ld$ and $\sF\lmt\sF^{\vee}$ are
quasi-inverse to each other. 
(See \cite[$\textrm{VII}_{\textrm{B}}$, 2.2.2.]{SGA3}) 

The Cartier-dual of a torus 
$\Trs\cong\left(\Gm\right)^{t}$
is a lattice of the same rank: 
$\Trsd\cong\Zint^{t}$, i.e.\ a torsion-free \'etale formal group 
(cf.\ \cite[(5.2)]{L}). 

Let $V$ be a finite dimensional $k$-vector space. 
The Cartier-dual of the vectorial group $\Vcl=\Spec\left(\Sym V^{\vee}\right)$
associated to $V$ is the infinitesimal formal group 
$\Vcld=\Spf\left(\widehat{\Sym V}\right)$, 
i.e.\ the completion w.r.t. $0$ of the vectorial
group associated to the dual $k$-vector space $V^{\vee}$ 
(cf.\ \cite[(5.2)]{L}).

\subsection{Duality of 1-Motives} 
\label{sub:Dual 1-Motive} 

The dual of an abelian variety $A$ is given by $\Ad=\Pic^0 A$. 
Unfortunately, 
there is no analogue duality theory for algebraic groups in general.
Instead, we embed the category of connected algebraic groups 
into the category of 1-motives 
by sending a connected algebraic group $G$ 
to the 1-motive $\left[0\ra G\right]$. 
The category of 1-motives admits a duality theory. 

\bThm 
\label{Struct} 
Let $L$ be a connected linear algebraic group, 
$A$ an abelian variety and $\Ad$ the dual abelian variety. 
There is a bijection 
\[ \Ext_{\Abk}\left(A,L\right)\,\simeq\,\Hom_{\Abk}\left(\Ld,\Ad\right) \;.
\] 
\eThm 

\bPf See \cite[Lemme (5.2.1)]{L} for a complete proof. \\ 
The bijection 
$\Phi:\Ext\left(A,L\right) \lra \Hom\left(\Ld,\Ad\right)$ 
is constructed as follows: 
Given an extension $G$ of $A$ by $L$. 
By Corollary \ref{fmlG determin} it suffices to determine 
the homomorphism $\Phi(G):\Ld\lra\Ad$ on the $k$-valued points 
and on the Lie-algebra of $\Ld$. 
Let $L=\Trs\tms\Vcl$ be the canonical splitting of $L$ 
into a direct product of a torus $\Trs$ and a vectorial group $\Vcl$ 
($\see$Theorem \ref{Thm linGrp}). Now 
\[ \begin{array}{c}
   \Ld(k)=\Trsd(k) %=\Homg(\Trs(k),\Gm(k))
   =\Homabk(\Trs,\Gm) \\
   \Lie\left(\Ld\right)=\Lie\left(\Vcld\right)=\Homk(\Lie(\Vcl),k)= 
   \Homabk(\Vcl,\Ga) 
   \end{array} 
\] 
i.e.\ $\chi\in\Ld(k)$ gives rise to a homomorphism $L\lra\Trs\lra\Gm$ and \\ 
$\chi\in\Lie\left(\Ld\right)$ to a homomorphism $L\lra\Vcl\lra\Ga$. 
Then the image of $\chi$ under $\Phi(G)$ is the push-out 
$\chi_* G \in \left\{ \begin{array}{l}
   \Ext\left(A,\Gm\right)=\Pic_A^0(k)=\Ad(k) \\ 
   \Ext\left(A,\Ga\right)=\Lie\left(\Pic_A^0\right)=\Lie\left(\Ad\right) 
                      \end{array} 
              \right.
$ 
\[ \xymatrix{  
   0\ar[r]&L\ar[d]_{\chi}\ar[r]&G\ar[d]\ar[r]&A\ar@{=}[d]\ar[r]&0 \\ 
   0 \ar[r] & \splG \ar[r] & \chi_* G \ar[r]  &  A \ar[r] & 0 \\ } 
\] 
where $\splG=\Gm$ or $\splG=\Ga$. 
\ePf 

\vspace{\vs} 
%\quad \\ 
A consequence of Theorem \ref{Struct} is the duality of 1-motives: \\ 
Let $M=\left[\sF\overset{\mu}\lra G\right]$ be a 1-motive,
and $0\ra L\ra G\ra A\ra0$ the canonical decomposition of $G$ 
($\see$Theorem \ref{Thm Chevalley}). 
By Theorem \ref{Struct} the composition $\bar{\mu}:\sF\lra G\lra A$
defines an extension $0\ra\Fd\ra G^{\bar{\mu}}\ra\Ad\ra0$
and the extension $G$ defines a homomorphism of sheaves of abelian 
groups $\bar{\mu}^{G}:\Ld\lra\Ad$. 
Then $\mu$ determines uniquely a factorization $\mu^G:\Ld\lra G^{\bar{\mu}}$ 
of $\bar{\mu}^{G}$, according to \cite[Proposition (5.2.2)]{L}, 
hence gives rise to 
\bDef The dual 1-motive of $M=\left[\sF\overset{\mu}\lra G\right]$ 
with $G\in\Ext\left(A,L\right)$ is the 1-motive 
$ M^{\vee}=\left[\Ld\overset{\mu^{G}}\lra G^{\bar{\mu}}\right]$. 
\eDef 

The double dual $M^{\vee\vee}$ of a 1-motive $M$ is canonically 
isomorphic to $M$ (see \cite[(5.2.4)]{L}).

\section{Universal Factorization Problem} 
\label{sec:Univ-Fact-Prbl}

Let $X$ be a projective variety over $k$ (an algebraically closed field 
of characteristic 0). 
The universal factorization problem may be outlined as follows: 
one is looking for a ``universal object'' $\sU$ and a rational map 
$u:X\dra\sU$ such that for every rational map 
$\phe:X\dra G$ to an algebraic group $G$ there is a 
unique homomorphism $h:\sU\lra G$ 
such that $\phe=h\circ u$ up to translation. 

The universal object $\sU$, if it exists, 
is not in general an algebraic group. 
For this aim a certain finiteness condition on the rational maps is needed. 
In this section we work out a criterion, 
for which categories $\Mr$ of rational maps from $X$ to algebraic groups 
one can find an algebraic group $\Alb_{\Mr}\left(X\right)$ 
satisfying the universal mapping property, 
and in this case we give a construction of $\Alb_{\Mr}\left(X\right)$.
The way of procedure was inspired by Serre's expos\'e \cite{S3}, 
where the case of semi-abelian varieties 
(see Examples \ref{Mr_D} and \ref{Serres genAlb}) is treated.

\subsection{Relative Cartier Divisors }

For the construction of universal objects as above we are concerned with 
the functor of families of Cartier divisors. 
This functor admits a natural transformation to the Picard functor, 
which describes families of line bundles, 
i.e.\ families of classes of Cartier divisors. 
References for \emph{effective} relative Cartier divisors are 
\cite[Lecture 10]{M2} and \cite[Section 8.2]{BLR}. 
Since the relative Cartier divisors we are concerned with 
are not necessarily effective, we give a short overview on this subject. 

\bDef \label{Lie} 
Let $F:\Algk\lra\Ab$ be a covariant functor, 
$R$ a $k$-algebra. 
Let $R_{\nil}:=k+\Nil(R)$ be the induced subring of $R$ 
and $\rho:R_{\nil}\lra k,\;\Nil\left(R\right)\ni n\lmt 0$
the augmentation. 
Define 
\begin{eqnarray*}
 \Inf\left(F\right)\left(R\right) & = & 
 \ker\Big(F\left(\rho\right):F\left(R_{\nil}\right)\lra F\left(k\right)\Big) \\ 
 \Lie\left(F\right) & = & \Inf\left(F\right)\left(\ke\right) \;.
\end{eqnarray*} 
\eDef 

\bNot \label{inn} 
If $F:\Algk\lra\Ab$ is a covariant functor, 
then we will use the expression $\alp\inn F$ in order to state that either 
$\alp\in F\left(k\right)$ or $\alp\in\Lie\left(F\right)$. 
\eNot

\subsubsection*{Functor of Relative Cartier Divisors }

Let $X$ and $Y$ be noetherian schemes of finite type over $k$. 
A \emph{Cartier divisor} on $X$ is by definition a global
section of the sheaf $\left.\sK_{X}^{*}\right/\sO_{X}^{*}$, where
$\sK_{X}$ is the sheaf of total quotient rings on $X$, and the star $^*$ 
denotes the unit groups. 
\[ \Div\left(X\right)=\Gamma\left(X,\left.\sK_{X}^{*}\right/\sO_{X}^{*}\right) 
\] 
is the \emph{group of Cartier divisors.} 

\bNot 
If $R$ is a $k$-algebra, 
the scheme $Y\times_{k}\Spec R$ is often denoted by $Y\tens R$. 
\eNot 

\bDef 
\label{relTotQuot}
Let $A$ be an $R$-algebra. 
The set 
\[ \Sy_{A/R} \,=\, 
   \left\{ f \in A \,\left| \begin{array}{c}
                        \textrm{$f$ not a zero divisor,}\\
                        A/fA \textrm{ is flat over R}\\
                        \end{array} 
   \right.\right\} 
\] 
is a multiplicative system in $A$. Then the localization of $A$
at $\Sy_{A/R}$ 
\[ \K_{A/R} \,=\, \Sy_{A/R}^{-1}A \] 
is called the \emph{total quotient ring of $A$} \emph{relative to $R$.} \\
Let $X\overset{\tau}\lra T$ be a scheme over $T$. The
sheaf $\sK_{X/T}$ associated to the presheaf %formed by the rings
\[  U \; \lmt \; 
    \K_{\sO_{X}\left(U\right)/(\tau^{-1}\sO_{T})\left(U\right)} \,=\,
    \Sy_{\sO_{X}\left(U\right)/(\tau^{-1}\sO_{T})\left(U\right)}^{-1}\sO_{X}\left(U\right) 
\] 
for open $U\subset X$, is called the \emph{sheaf of total quotient
rings of $X$ relative to $T$.} 
\eDef 

\bRmk 
Let $A$ be a finitely generated flat $R$-algebra, 
where $R$ is a noetherian ring, 
$f\in A$ a non zero divisor. 
Notice that $A/fA$ is flat over $R$ if and only if for all homomorphisms 
$R\lra S$ the image of $f$ in $A\tens_R S$ is a non zero divisor. 
Equivalently, $f$ is not contained in any associated prime ideal 
of $A\tens_R k(p)$ for all $p\in\Spec R$ 
(cf.\ \cite[Lecture 10]{M2}). 
\eRmk 

\bPrp 
\label{Div_Y}
For a $k$-algebra $R$ let 
\[ \Divf_{Y}\left(R\right) =
   \Gamma\left(Y\tens R,\left.\sK_{Y\tens R/R}^*\right/\sO_{Y\tens R}^*\right) \;. 
\]
Then the assignment \quad $R\lmt\Divf_{Y}\left(R\right)$ \; 
defines a covariant functor 
\[ \Divf_{Y}: \Algk \lra \Ab 
\] 
from the category of $k$-algebras to the category of abelian groups. 
\ePrp 

\bRmk 
\label{Div_Y criterium}
For each $k$-algebra $R$ we have 
\[ \Divf_{Y}\left(R\right)=\left\{ 
   \begin{array}{c}
     \textrm{Cartier divisors $\sD$ on $Y\times_{k}\Spec R$}\\
     \textrm{which define Cartier divisors 
             $\sD_{p}$ on $Y\times\left\{ p\right\} $} \\
     \forall p\in\Spec R 
   \end{array} 
\right\} 
\] 
and for a homomorphism $h:R\lra S$ in $\Algk$ the induced
homomorphism $\Divf_{Y}\left(h\right):
\Divf_{Y}\left(R\right)\lra\Divf_{Y}\left(S\right)$
in $\Ab$ is the pull-back of Cartier divisors on $Y\times_{k}\Spec R$
to those on $Y\times_{k}\Spec S$. 
\eRmk 

\bPrp \label{Divf(k), Lie(Divf)} 
Let $R$ be a $k$-algebra, 
$R_{\nil}=k+\Nil\left(R\right)$ 
the induced ring. 
There is a canonical isomorphism of abelian groups 
\[ \Inf\left(\Divf_{Y}\right)\left(R\right) = 
   \Lie\left(\Divf_{Y}\right)\tens_{k}\Nil\left(R\right) 
\] 
where 
$\Lie\left(\Divf_{Y}\right) = \Gamma\left(Y,\left.\sK_{Y}\right/\sO_{Y}\right)$. 
\ePrp 

\bPf  $R_{\nil}=k\left[\Nil\left(R\right)\right]$ is a local %Artinian
ring with only prime ideal $\Nil\left(R\right)\in\Spec R_{\nil}$,
and it consists of zero divisors only. Therefore it holds 
$\sK_{Y\tens R_{\nil}/R_{\nil}}^{*}=\sK_{Y\tens R_{\nil}}^{*}
 =\sK_{Y}\left[\Nil\left(R\right)\right]^{*}=
\sK_{Y}^{*}+\sK_{Y}\tens_{k}\Nil\left(R\right)$. 
Hence we obtain an exact sequence 
\[ 1\lra
\frac{\sO_Y^*+\sK_{Y}\tens_{k}\Nil\left(R\right)}
{\sO_Y^*+\sO_{Y}\tens_{k}\Nil\left(R\right)}
\lra\frac{\sK_{Y\tens R_{\nil}}^{*}}{\sO_{Y\tens R_{\nil}}^{*}}
\lra\frac{\sK_{Y}^{*}}{\sO_{Y}^{*}}\lra1 \;. 
\] 
Now we have canonical isomorphisms 
\[
\frac{\sO_Y^*+\sK_{Y}\tens_{k}\Nil\left(R\right)}
{\sO_Y^*+\sO_{Y}\tens_{k}\Nil\left(R\right)} \cong 
\frac{1+\sK_{Y}\tens_{k}\Nil\left(R\right)}{1+\sO_{Y}\tens_{k}\Nil\left(R\right)}
\cong\frac{\sK_{Y}\tens_{k}\Nil\left(R\right)}{\sO_{Y}\tens_{k}\Nil\left(R\right)}
\cong\frac{\sK_{Y}}{\sO_{Y}}\tens_{k}\Nil\left(R\right) 
\] 
where the second isomorphism is given by $\exp^{-1}$. 
Applying the global section functor $\Gamma\left(Y,\_\right)$ yields \[
0\lra\Gamma\left(\frac{\sK_{Y}}{\sO_{Y}}\tens_{k}\Nil\left(R\right)\right)\lra
\Gamma\left(\frac{\sK_{Y\tens R_{\nil}}^{*}}{\sO_{Y\tens R_{\nil}}^{*}}\right)\lra
\Gamma\left(\frac{\sK_{Y}^{*}}{\sO_{Y}^{*}}\right) \;. 
\] 
Here 
$\Gamma\left(\left.\sK_{Y\tens R_{\nil}}^{*}\right/\sO_{Y\tens R_{\nil}}^{*}\right) = 
 \Divf_{Y}\left(R_{\nil}\right)$
and 
$\Gamma\left(\left.\sK_{Y}^{*}\right/\sO_{Y}^{*}\right) = 
 \Divf_{Y}\left(k\right)$, 
therefore \[
\Inf\left(\Divf_{Y}\right)\left(R\right)=
\Gamma\left(\left.\sK_{Y}\right/\sO_{Y}\right)\tens_{k}\Nil\left(R\right) \;. 
\] 
In particular for $R=\ke$, as $\Nil\left(\ke\right)=\varepsilon k\cong k$, 
we have \[
\Lie\left(\Divf_{Y}\right)=\Inf\left(\Divf_{Y}\right)\left(\ke\right)\cong
\Gamma\left(\left.\sK_{Y}\right/\sO_{Y}\right) 
\] 
and hence \; 
$\Inf\left(\Divf_{Y}\right)\left(R\right)\cong
\Lie\left(\Divf_{Y}\right)\tens_{k}\Nil\left(R\right)$. 
\ePf 

\bDef \label{Support} 
If $D\in\Divf_{Y}\left(k\right)$, then 
$\Supp\left(D\right)$ denotes the locus of zeros and poles of local sections 
$\left(f_{\alpha}\right)_{\alpha}$ of $\sK_{Y}^{*}$ representing 
$D \in \Gamma\left(\left.\sK_{Y}^{*}\right/\sO_{Y}^{*}\right)$. \\
If $\delta\in\Lie\left(\Divf_{Y}\right)$, then $\Supp\left(\delta\right)$
denotes the locus of poles of local sections 
$\left(g_{\alpha}\right)_{\alpha}$ of $\sK_{Y}$ 
representing $\delta\in\Gamma\left(\left.\sK_{Y}\right/\sO_{Y}\right)$. 
\eDef 

\bDef \label{Supp(F)} 
Let $\fmlG$ be a subfunctor of $\Divf_Y^0$ which is a formal group. 
Then $\Supp(\fmlG)$ is defined to be the union of $\Supp(D)$ 
for $D\in\fmlG(k)$ and $D\in\Lie(\fmlG)$. 
\eDef  
$\Supp \lp \fmlG \rp$ is a closed subscheme of codimension 1 in $Y$, since
$\fmlG\left(k\right)$ and $\Lie\left(\fmlG\right)$ are both finitely
generated. 

\bDef \label{Dec_Y,V} 
For a morphism $\sigma: V \lra Y$ of 
varieties define $\Decf_{Y,V}$ to be the subfunctor of $\Divf_{Y}$ 
consisting of families of those Cartier divisors 
that do not contain any component of the image $\sigma(V)$. 
\eDef  

\bPrp \label{_.V} Let $\sigma: V \lra Y$ be a morphism of 
varieties. Then the pull-back of Cartier divisors $\sigma^*$ 
induces a natural transformation of functors 
\[ \_\cut V:\Decf_{Y,V}\lra\Divf_{V} \;. 
\] 
\ePrp 

\bPrp \label{fmlG->Div} 
Let $\fmlG$ be a formal group. Then each pair $(a,l)$ of 
a homomorphism of abelian groups $a:\fmlG(k)\lra\Divf_Y(k)$ and 
a $k$-linear map $l:\Lie\left(\fmlG\right)\lra\Lie\left(\Divf_Y\right)$ 
determines uniquely a natural transformation $\trafo:\fmlG\lra\Divf_Y$ 
with $\trafo(k)=a$ and $\Lie(\trafo)=l$. 
Moreover, the image of $\fmlG$ in $\Divf_Y$ is also a formal group. 
\ePrp 

\bPf We construct a natural transformation with the required property 
by giving homomorphisms $\trafo(R):\fmlG(R)\lra\Divf_Y(R)$, $R\in\Algk$. 
As $\fmlG$ and $\Divf_Y$ both commute with direct products, 
we may reduce to $k$-algebras $R$ with $\Spec R$ connected. 
In this case 
$\fmlG(R)=\fmlG(k)\tms\big(\Lie\left(\fmlG\right)\tens_k\Nil(R)\big)$, 
cf.\ Theorem \ref{Str-fml-Grp} and Remark \ref{G_inf(R)}. 
Then $(a,l)$ yields a homomorphism 
\[ h: \fmlG(k)\tms\big(\Lie\left(\fmlG\right)\tens_k\Nil(R)\big) \lra 
       \Divf_Y(k)\tms\big(\Lie\left(\Divf_Y\right)\tens_k\Nil(R)\big) \;. 
\] 
We have functoriality maps $\Divf_Y(k)\ra\Divf_Y(R)$ and 
$\Divf_Y\left(R_{\nil}\right) \ra \Divf_Y(R)$, and 
$\Lie\left(\Divf_Y\right)\tens_k\Nil(R) \cong \Inf\left(\Divf_Y\right)(R) 
  \subset \Divf_Y\left(R_{\nil}\right)$. 
Then composition with $h$ gives $\trafo(R)$. 
It is straightforward that the homomorphisms 
$\trafo(R)$ yield a natural transformation. 
The definition of $\Inf\left(\_\right)$ implies that 
$\Divf_Y(k)\tms\Inf\left(\Divf_Y\right)(R) \lra \Divf_Y(R)$ is injective. 
Then the second assertion follows from the structure of formal groups, 
cf.\ Corollary \ref{fmlG determin}. 
\ePf

\subsubsection*{Picard Functor }

Although the Picard functor is an established object in algebraic geometry, 
we give a summary of facts that we need in the following. 
References for the Picard functor are e.g.\ \cite[Chapter 8]{BLR}, 
\cite[Lectures 19-21]{M2} or \cite{K}. 

Let $Y$ be a projective scheme over $k$. 
The isomorphism classes of line bundles on a $k$-scheme $X$ form a group
$\Pic\left(X\right)$, the \emph{(absolute) Picard group} of $X$,
which is given by 
$\Pic\left(X\right)=\H^{1}\left(X,\sO_{X}^{*}\right)$.  
The \emph{(relative) Picard functor} $\Picf_Y$ 
from the category of $k$-algebras 
to the category of abelian groups 
is defined by 
\[ \Picf_{Y}\left(R\right)=
   \left.\Pic\left(Y\times_{k}\Spec R\right)\right/\Pic\left(\Spec R\right) 
\] 
for each $k$-algebra $R$. 
In other words, $\Picf_Y(R)$ is the abelian group of line bundles on 
$Y\tens R$ modulo line bundles that arise as a pull-back of a line bundle 
on $\Spec R$. 

Similarly to the proof of Proposition \ref{Divf(k), Lie(Divf)} 
one shows 

\bPrp  \label{Picf(k), Lie(Picf)}
Let $R$ be a $k$-algebra, $R_{\nil}=k+\Nil\left(R\right)$ the induced ring. 
There is a canonical isomorphism of abelian groups  
\[ \Inf\left(\Picf_{Y}\right)\left(R\right) = 
   \Lie\left(\Picf_{Y}\right)\tens_{k}\Nil\left(R\right) 
\] 
where $\Lie\left(\Picf_{Y}\right) = \H^{1}\left(Y,\sO_{Y}\right)$. 
\ePrp 

Since a scheme over an algebraically closed field $k$ admits a section, 
the fppf-sheaf associated to the Picard functor on $Y$ coincides 
with the relative Picard functor $\Picf_Y$, 
see \cite[Section 8.1 Proposition 4]{BLR}. 

\bThm \label{Pic_scheme} 
Let $Y$ be an integral projective $k$-scheme of finite type. 
Then the Picard functor $\Picf_{Y}$ is represented by
a $k$-group scheme $\Pic_{Y}$, which is called the 
\emph{Picard scheme of $Y$}. 
\eThm 

\bPf \cite[No.~232, Theorem 2]{FGA} or  
\cite[Section 8.2, Theorem 1]{BLR} or 
\cite[Theorem 4.8, Theorem 4.18.1]{K}. 
\ePf 

\bDef Let $M,N$ be line bundles on $Y$. Then $M$ is said to be
\emph{algebraically equivalent} to $N$, %(in signs $M\sim N$), 
if there exists a connected $k$-scheme $C$, 
a line bundle $\sL$ on $Y\times_{k}C$ and closed points $p,q\in C$ 
such that $\sL|_{Y\tms\left\{ p\right\} }=M$ and $\sL|_{Y\tms\left\{ q\right\} }=N$. 
\eDef 

If $\Picf_{Y}$ is represented by a scheme $\Pic_{Y}$, 
then a line bundle $L$ on $Y$ is algebraically equivalent to $\sO_Y$ 
if and only if $L$ lies in the connected component $\Pic_{Y}^{0}$ 
of the identity of $\Pic_{Y}$. 
Then $\Pic_Y^0$ represents the functor $\Picf_Y^0$ 
which assigns to $R\in\Algk$ 
the abelian group of line bundles $\sL$ on $Y\tens R$ 
with $\sL|_{Y\times\left\{ t\right\} }$ algebraically equivalent to $\sO_{Y}$
for each $t\in\Spec R$, 
modulo line bundles that come from $\Spec R$. 

\bThm 
\label{proj Pic^0} 
Let $Y$ be an integral projective $k$-scheme of finite type. 
Then $\Picf_{Y}^{0}$ is represented by a quasi-projective 
$k$-group scheme $\Pic_{Y}^{0}$. 
If $Y$ is also normal, then $\Pic_{Y}^{0}$ is projective. 
\eThm 

\bPf  Follows from Theorem \ref{Pic_scheme} 
and \cite[Theorem 5.4]{K}. 
\ePf 

\vspace{\vs} 

A normal projective variety $Y$ over $k$ is the disjoint union of its 
irreducible components %$Z \in \Cp(Y)$ 
(see \cite[Chapter 3, \S 9, Remark p.~64]{Mm}). 
Applying Theorem \ref{proj Pic^0} to each irreducible component $Z$ of $Y$ 
yields that $\Picf_Y^0$ is represented by the product of 
the connected projective $k$-group schemes $\Pic_Z^0$. 
In characteristic 0 a connected projective $k$-group scheme 
is an abelian variety. 
We obtain 

\bCor \label{normal Pic^0} Let $Y$ be a normal projective variety over $k$. 
Then $\Picf_Y^0$ is represented by an abelian variety $\Pic_Y^0$. 
\eCor

\subsubsection*{Transformation $\Divf_{Y}\lra\Picf_{Y}$ }

Let $Y$ be a $k$-scheme and $R$ be a $k$-algebra. 
Consider the exact sequence $\Seq\left(R\right)$: 
\[ 1\lra\sO_{Y\tens R}^{*}\lra\sK_{Y\tens R}^{*}
   \lra\left.\sK_{Y\tens R}^{*}\right/\sO_{Y\tens R}^{*}\lra1 \;. 
\] 
In the corresponding long exact sequence $\H^{\bullet}\Seq\left(R\right)$:
\[ \lra\H^{0}\left(\sK_{Y\tens R}^{*}\right)\lra
   \H^{0}\left(\left.\sK_{Y\tens R}^{*}\right/\sO_{Y\tens R}^{*}\right)
   \lra\H^{1}\left(\sO_{Y\tens R}^{*}\right)\lra
   \H^{1}\left(\sK_{Y\tens R}^{*}\right)\lra 
\] 
the connecting homomorphism 
$\delta^{0}\left(R\right):
 \H^{0}\left(\left.\sK_{Y\tens R}^{*}\right/\sO_{Y\tens R}^{*}\right)
 \lra\H^{1}\left(\sO_{Y\tens R}^{*}\right)$ 
gives a natural transformation 
$\Div\left(Y\tens R\right)\lra\Pic\left(Y\tens R\right)$.
Composing this transformation with the injection 
$\Divf_{Y}\left(R\right)\inj\Div\left(Y\tens R\right)$
and the projection $\Pic\left(Y\tens R\right)\sur\Picf_{Y}\left(R\right)$
yields a natural transformation 
\[ \cl:\Divf_{Y}\lra\Picf_{Y} \;. 
\] 

\bDef \label{Div^0} 
Let $\Divf_{Y}^{0}$
be the subfunctor of $\Divf_{Y}$ defined by 
\;$ \Divf_{Y}^{0}\left(R\right)=\cl^{-1}\left(\Picf_{Y}^{0}\left(R\right)\right) 
$\; 
for each $k$-algebra $R$. 
\eDef

\subsection{Categories of Rational Maps to Algebraic Groups} 
\label{sub:Categories-of-Rational} 

Let $Y$ be a regular projective variety over $k$ 
(an algebraically closed field of characteristic 0). 
Algebraic groups are always assumed to be connected, unless stated otherwise. 

\bNot  $\splG$ stands for one of the groups $\Gm$ or $\Ga$. 
\eNot 

\bLem \label{div(s)}Let $P$ be a principal $\splG$-bundle over $Y$.
Then a local section $\sigma:U\subset Y\lra P$ determines
uniquely a divisor $\dv_{\Gm}\left(\sigma\right)\in\Divf_{Y}(k)$ on $Y$, 
if $\splG=\Gm$, or an infinitesimal deformation 
$\dv_{\Ga}\left(\sigma\right)\in\Lie\left(\Divf_{Y}\right)$
of the zero divisor on $Y$, if $\splG=\Ga$. 
\eLem 

\bPf For 
$V = \left\{ \begin{array}{cc}
             k & \textrm{if } \splG=\Gm \\ 
             \ke & \textrm{if } \splG=\Ga \\ 
             \end{array}
     \right.$ 
let $\lambda: \splG(k) \lra \Gl\left(V\right)$ be the representation 
of $\splG$ given by 
$l \lmt \left\{ \begin{array}{cc}
                       l & \textrm{if }\splG=\Gm\\
                       1+\varepsilon l & \textrm{if }\splG=\Ga
                       \end{array} 
               \right.$.
Let $\sV=P\amalg_{\left[\splG,\lambda\right]}V$ be the vector-bundle associated to 
$P$ of fibre-type $V$. Denote by $s\in\Gamma\left(U,\sV\right)$ the image 
of $\sigma$ under the map $P \lra \sV$ 
induced by $\lambda$ on the fibres. 
There is an effective divisor $H$, supported on $Y\setminus U$,
such that the local section $s\in\Gamma\left(U,\sV\right)$ extends to
a global section of $\sV\left(H\right)=\sV\tens_{\sO}\sO_Y(H)$.
Local trivializations $\sV|_{U_{\alpha}} \overset{\sim}\lra\sO_{U_{\alpha}}\tens_{k}V$ 
induce local isomorphisms  
$\Phi_{\alpha}:\sV(H)|_{U_{\alpha}} \overset{\sim}\lra\sO(H)|_{U_{\alpha}}\tens_{k}V$ 
of the twisted bundles. 
Then the local sections 
$\Phi_{\alpha}\left(s\right)\in\Gamma\left(U_{\alpha},\sO_{Y}(H)\tens_{k}V\right)$
yield a divisor $\dv_{\Gm}\left(\sigma\right)$ on $Y$ if $V=k$, 
and a divisor $\dv_{\Ga}\left(\sigma\right)$ 
on $Y\left[\varepsilon\right]=Y\times_{k}\Spec\ke$ if $V=\ke$. 
The description of $\Divf_Y$ from Remark \ref{Div_Y criterium} shows that 
$\dv_{\splG}\left(\sigma\right)\in\Divf_Y(V)$. 
In the second case, the image of $\lambda:\Ga(k)\lra\Gl(\ke)$ 
lies in $1+\eps k$. 
This implies that the restriction of $\dv_{\Ga}\left(\sigma\right)$ to $Y$ 
is the zero divisor, 
thus $\dv_{\Ga}\left(\sigma\right) \in \Lie\left(\Divf_{Y}\right)$. 
\ePf 

\vspace{\vs} 

Let $\phe:Y \dra G$ be a rational map to an algebraic
group $G$ with canonical decomposition 
$0\ra L\ra G\overset{\rho}\ra A\ra0$.
Since a rational map to an abelian variety is defined at every regular point 
(see \cite[Chapter II, \S 1, Theorem 2]{La}), 
the composition $Y \overset{\phe} \dra G \overset{\rho} \lra A$
extends to a morphism $\overline{\phe}:Y \lra A$. 
Let $G_{Y}=G \times_{A} Y$ be the fibre-product of $G$ and $Y$ 
over $A$. The graph-morphism $\phe_{Y}:U\subset Y \lra 
G \times_{A} Y,\quad y\lmt\left(\phe\left(y\right),y\right)$
of $\phe$ is a section of the $L$-bundle $G_{Y}$ over $Y$. 
Each $\lambda\inn\Ld$ (see Notation \ref{inn}, where we consider
$\Ld$ as a functor on $\Algk$) gives rise to a homomorphism 
$\lambda:L\lra\splG$ (cf.\ proof of Theorem \ref{Struct}). 
Then the composition of $\phe_{Y}$ 
with the push-out of $G_{Y}$ via $\lambda$ gives a local section 
$\phe_{Y,\lambda}:U\subset Y\lra\lambda_{*}G_{Y}$
of the $\splG$-bundle $\lambda_{*}G_{Y}$ over $Y$: 
\[ \xymatrix{ 
   \splG \ar[dd] && L \ar[ll]_{\lambda} \ar[dd] \ar@{=}[rr] && L \ar[dd] \\  
   && && \\  
   \lambda_{*} G_Y \ar@<-0.7ex>[dd] && 
   G_Y \ar[ll] \ar@<-0.7ex>[dd] \ar[rr]  &&  G \ar[dd]^{\rho} \\  
   && && \\  
   Y \ar@{-->}@<-0.7ex>[uu]_{\phe_{Y,\lambda}} \ar@{=}[rr] && 
   Y \ar@{-->}@<-0.7ex>[uu]_{\phe_Y} \ar@{-->}[uurr]_{\phe} \ar[rr] && A } 
\] 
Lemma \ref{div(s)} says that the local section $\phe_{Y,\lambda}$ 
determines a unique divisor or deformation 
$\dv_{\splG}\left(\phe_{Y,\lambda}\right)\in\left\{ 
\begin{array}{cc}
\Gamma\left(\left.\sK_{Y}^{*}\right/\sO_{Y}^{*}\right) & \textrm{if }\splG=\Gm\\
\Gamma\left(\left.\sK_{Y\left[\varepsilon\right]}^{*}\right/\sO_{Y\left[\varepsilon\right]}^{*}
 \right) & \textrm{if }\splG=\Ga\end{array}\right.$. 
Now the bundle $\lambda_{*}G_{Y}$ comes from an extension of algebraic groups, 
and $\Extfabk(A,\Gm) \cong \Picf^0_A$, 
hence it is an element of $\Picf^0_Y(k)$ if $\splG=\Gm$, 
or of $\Lie\left(\Picf^0_Y\right)$ if $\splG=\Ga$. 
Therefore $\dv_{\splG}\left(\phe_{Y,\lambda}\right)$ is a divisor in 
$\Divf_{Y}^{0}(k)$ if $\lambda\in L\left(k\right)$,
or an element of $\Lie\left(\Divf_{Y}^{0}\right)$ 
if $\lambda\in\Lie\left(L\right)$. 

\bPrp \label{induced Trafo} 
Let $G\in\Ext\left(A,L\right)$ be an
algebraic group and $\phe: Y \dra G$ a rational map.
Then $\phe$ induces a natural transformation of functors 
$\trafo_{\phe}: \Ld \lra \Divf_{Y}^{0}$.
\ePrp  

\bPf  
The construction above yields a homomorphism of abelian groups
$\Ld\left(k\right)\ra\Divf_{Y}^{0}\left(k\right)$, 
$\lambda\mpt\dv_{\Gm}\left(\phe_{Y,\lambda}\right)$
and a $k$-linear map 
$\Lie\left(L\right)\ra\Lie\left(\Divf_{Y}^{0}\right),$
$\lambda\mpt\dv_{\Ga}\left(\phe_{Y,\lambda}\right)$. 
Then by Proposition \ref{fmlG->Div} 
this extends to a natural transformation $\Ld \lra \Divf_{Y}^{0}$. 
\ePf 

\bDef \label{CatMr} 
A category $\Mr$ is called a 
\emph{category of rational maps from $Y$
to algebraic groups}, if objects and morphisms of $\Mr$ 
satisfy the following conditions: 
The objects of $\Mr$ are rational maps $\phe: Y \dra G$, 
where $G$ is an algebraic group, 
such that $\phe(U)$ generates a \emph{connected} 
algebraic subgroup of $G$ for any open set $U \subset Y$ 
on which $\phe$ is defined. 
The morphisms of $\Mr$ between two objects $\phe: Y \dra G$
and $\psi: Y \dra H$ are given by the set of all homomorphisms of 
algebraic groups $h:G \lra H$ such that $h \circ \phe = \psi$, 
i.e.\ the following diagram commutes: 
\[  \xymatrix{  & Y \ar@{-->}[dl]_{\phe} \ar@{-->}[dr]^{\psi} &  \\  
                G \ar[rr]^{h} & & H \\ 
             } 
\] 
\eDef 

\bRmk \label{EquCatMr} Let $\phe:Y\dra G$ and $\psi:Y\dra H$ 
be two rational maps from $Y$ to algebraic groups. Then Definition 
\ref{CatMr} implies that for any category $\Mr$ of rational maps 
from $Y$ to algebraic groups containing $\phe$ and $\psi$ as objects 
the set of morphisms $\Hom_{\Mr}(\phe,\psi)$ is the same.  
Therefore two categories $\Mr$ and $\Mr'$ of rational maps from $Y$ to 
algebraic groups are equivalent if every object of $\Mr$ is isomorphic to 
an object of $\Mr'$. 
\eRmk 

\bDef \label{H_I} We denote by $\Insg \Mr$ 
the localization of a category $\Mr$ 
of rational maps from $Y$ to algebraic groups 
at the system of injective homomorphisms 
among the morphisms of $\Mr$. 
\eDef 

\bRmk In $\Insg \Mr$ we may assume  
for an object $\phe:Y\dra G$ that $G$ is generated by $\phe$: 
The inclusion $\langle\im\phe\rangle \inj G$ of 
the subgroup $\langle\im\phe\rangle$ generated by $\phe$ 
is an injective homomorphism, hence $\phe:Y\dra G$ is isomorphic to 
$\phe:Y\dra \langle\im\phe\rangle$, 
if $(\phe:Y\dra \langle\im\phe\rangle) \in \Mr$. 
\eRmk

\bDef \label{Mav} The category of rational maps from $Y$ 
to abelian varieties is denoted by $\Mav$. 
\eDef 

\bRmk The objects of $\Mav$ are in fact \emph{morphisms} 
from $Y$, since a rational map from a regular variety $Y$ to an abelian 
variety $A$ extends to a morphism from $Y$ to $A$. 
\eRmk 

\bDef \label{Mr_sW} Let $\fmlG$ be a subfunctor of $\Divf_{Y}^{0}$ 
which is a formal group. Then $\Mr_{\fmlG}$ denotes the category of 
those rational maps $\phe:Y\dra G$ from $Y$ to algebraic groups 
for which the images of the natural transformations 
$\trafo_{\phe}:\Ld\lra\Divf_{Y}^{0}$
($\see$Proposition \ref{induced Trafo}) lie in $\fmlG$, 
i.e.\ which induce homomorphisms of formal groups $\Ld\lra\fmlG$,
where $0\ra L\ra G\ra A\ra0$ is the canonical decomposition of $G$. 
\[ \Mr_{\fmlG} = \{ \phe:Y\dra G \;|\; \im \trafo_{\phe} \subset \fmlG \}
\]
\eDef 

\bExm \label{Mr_0} For the category $\Mr_{0}$ associated to the trivial
formal group $\{0\}$ the localization $\Insg \Mr_0$ ($\see$Definition 
\ref{H_I}) is equivalent to the localization $\Insg \Mav$ of the category 
of morphisms from $Y$ to abelian varieties: 

Let $\phe:Y\dra G$ be a rational map to an algebraic group $G\in\Ext(A,L)$, 
and assume $\phe$ generates a connected subgroup of $G$. 
Then the following conditions are equivalent: \\ 
\begin{tabular}{rl} 
(i) & The transformation $\trafo_{\phe}:\Ld\lra\Divf_Y^0$ induced by $\phe$  
      has image $\{0\}$. \\
(ii) & The section $\phe_{Y}$ to the $L$-bundle $G_{Y}$ over $Y$ 
      extends to a global section. \\
\end{tabular} \\ 
In this case $\phe$ is defined on the whole of $Y$, 
i.e.\ is a morphism from $Y$ to $G$. 
Since $Y$ is complete, $\im\phe$ is a complete subvariety of $G$. 
Then the subgroup $\langle\im\phe\rangle$ of $G$ is also complete 
(cf.\ \cite[Lemma 1.10 (ii)]{ESV}) and connected, hence an abelian variety. 
\eExm 

\bDef \label{Resolution of Sing} 
Let $X$ be a (singular) projective variety. 
A morphism of varieties $\pi:\Xt\lra X$ is called a 
\emph{resolution of singularities} for $X$, 
if $\Xt$ is nonsingular and 
$\pi$ is a proper birational morphism 
which is an isomorphism over the nonsingular points of $X$. 
\eDef 

\bExm \label{M_0^X} 
Let $X$ be a singular projective variety and $Y=\Xt$, 
where $\pi:\Xt\lra X$ is a projective resolution of singularities. 
Denote by $\Mro{X}$ the category of rational maps $\phe:X\dra G$ 
whose associated map on 0-cycles of degree 0 $\ZOo{U}\lra G(k)$ 
(where $U$ is the open set on which $\phe$ is defined) 
factors through the homological Chow group of 0-cycles of degree 0 
modulo rational equivalence $\A_0(X)^0$. 
The definition of $\A_0(X)$ (see \cite[Section 1.3, p.10]{F}) 
implies that such a rational map is necessarily defined on the 
whole of $X$, i.e.\ is a morphism $\phe:X\lra G$. 
Then the composition $\phe\circ\pi:Y\lra G$ is an object of $\Mr_0$. 
Thus $\Mro{X}$ is a subcategory of the category $\Mr_0$ 
from Example \ref{Mr_0}. 
\eExm 

\bExm \label{Mr_D} Let $D$ be a reduced effective divisor on $Y$. 
Let $\Mr^D$ be the category of rational maps from $Y$ to 
semi-abelian varieties (i.e.\ extensions of an abelian variety by a torus) 
which are regular away from $D$. 

Let $\fmlG_{D}$ be the formal group 
whose \'etale part is given by divisors in $\Divf_{Y}^{0}\left(k\right)$ 
with support in $\Supp\left(D\right)$ and 
whose infinitesimal part is trivial. 

Then $\Mr^D$ is equivalent to $\Mr_{\fmlG_{D}}$: 
For a rational map $\phe:Y\dra G$ the induced sections
$\phe_{Y,\lambda}$ determine divisors in $\Divf_{Y}^{0}\left(k\right)$
supported on $\Supp\left(D\right)$ for all $\lambda\inn\Ld$ 
(where $L$ is the largest linear subgroup of $G$) 
if and only if $L$ is a torus, i.e.\ it consists of several copies 
of $\Gm$ only, and $\phe$ is regular on $Y\setminus\Supp\left(D\right)$.
\eExm 

\bExm \label{Mr_del} Let $Y=C$ be a smooth projective curve, 
$\mdl=\sum_{i}n_{i}\, p_{i}$ with $p_{i}\in C$, $n_{i}$ integers $\geq1$, 
an effective divisor on $C$ 
and let $\val_{p}$ be the valuation attached to the point $p\in C$.

Let $\Mr^{\mdl}$ be the category of those rational maps $\phe:Y\dra G$ 
such that for all $f\in\sK_{C}$ it holds: 
\[ \val_{p_{i}}\left(1-f\right)\geq n_{i}\quad\forall i\quad
   \Lra\quad\phe\left(\dv\left(f\right)\right)=0 \;. 
\] 
Let $\fmlG_{\mdl}$ be the formal group defined by 
\begin{eqnarray*} 
   \fmlG_{\mdl}\left(k\right) & = & 
   \left\{ \left.\sum_{i}l_{i}\, p_{i}\;\right|\;l_i\in\Zint,\;
   \sum_{i}l_{i}=0 \right\} \\ 
   \Lie \left( \fmlG_{\mdl} \right) & = & 
   \Gam\left(\left.\sO_C\left(\sum_i(n_i-1)\,p_i\right)
             \right/\sO_C\right) \;. 
\end{eqnarray*} 
Then $\Mr^{\mdl}$ is equivalent to $\Mr_{\fmlG_{\mdl}}$. 
By constructing a singular curve associated to the modulus $\mdl$ 
(see \cite[Chapter IV, No.~4]{S}), this turns out to be a special 
case of Example \ref{Mr_Y/X} 
(cf.\ Lemma \ref{IDiv_Y/X^0 repr} for the computation of 
$\Lie \left( \fmlG_{\mdl} \right)$). 
\eExm 

\bExm \label{Mr_Omega} 
Let $Y$ be a smooth projective variety over $\Cplx$, $D$ a divisor on $Y$ 
with normal crossings, and let $U=Y\setminus\Supp(D)$. 
Let $W\subset\Gam\left(U,\Oma_U^1\right)^{\der=0}$ be a finite dimensional 
$k$-vector space containing $\Gam\left(Y,\Oma_Y^1[\log D]\right)$. 
Let $\Mr^W$ be the category of those morphisms 
$\phe:U\lra G$ from $U$ to algebraic groups for which 
$ \phe^* \left( \Lie\left(G\right)^{\vee} \right) \subset W $, 
where $\Lie\left(G\right)^{\vee} = \Gam\left(G,\Oma_G^1\right)^{\cst}$ 
is the $k$-vector space of translation invariant regular 1-forms on $G$. 
A 1-form $\oma\in W$ determines a deformation of the zero divisor 
$\del(\oma)\in\Gam\left(\left.\sK_Y\right/\sO_Y\right)$ as follows: 
There are a covering $\{V_i\}_i$ of $Y$ and regular functions 
$f_i\in\sO_Y(U\cap V_i)$ such that $\oma-\der f_i$ is regular on $V_i$ 
(see \cite[VI.4 proof of Lemma 7]{FW}). 
Define $\del(\oma)=\left[\left(f_i\right)_i\right]\in 
 \Gam\left(\left.\sK_Y\right/\sO_Y\right)$. 

Let $\fmlG_W$ be the formal group determined by 
\begin{eqnarray*} 
   \fmlG_W\left(k\right) & = & 
      \left\{ \left. D' \in \Divf_Y^0 \left(k\right) \, \right| \, 
      \Supp(D')  \subset \Supp(D) \right\} \\ 
   \Lie \left( \fmlG_W \right) & = & 
      \im\big(\del:W\lra\Gam\left(\left.\sK_Y\right/\sO_Y\right)\big) \;. 
\end{eqnarray*} 
Then $\Mr^W$ is equivalent to $\Mr_{\fmlG_W}$. 
This follows from the construction of the generalized Albanese variety 
of Faltings and W\"ustholz (see \cite[VI.2. Satz 6]{FW}). 
\eExm 

\bExm \label{Mr_Y/X} Let $X$ be a singular projective variety and
$Y=\Xt$, where $\pi:\Xt\lra X$ is a projective resolution of singularities.
A rational map $\phe:X\dra G$ which is regular on the
regular locus $X_{\reg}$ of $X$ can also be considered as a rational
map from $Y$ to $G$. 
Let $\Mr^{\CHOo{X}}$ be the category of morphisms $\phe:X_{\reg}\lra G$
which factor through a homomorphism of groups 
$\CHOo{X}\lra G(k)$, see Definition \ref{MrCH(X)} 
(cf.\ \cite[Definition 1.14]{ESV} for the notion of \emph{regular 
homomorphism}). 
Let $\Divf_{\Xt/X}^{0}$ be the formal group given by the kernel 
of the push-forward $\pi_*$ ($\see$Proposition \ref{Div_Y/X^0}). \\ 
Then $\Mr^{\CHOo{X}}$ is equivalent to $\Mr_{\Divf_{\Xt/X}^{0}}$. 
This is the subject of Section \ref{sec:Rat-Maps CH_0(X)_deg0}. 
\eExm

\subsection{Universal Objects} 
\label{sub:Universal-Objects}

Let $Y$ be a regular projective variety over $k$ 
(an algebraically closed field of characteristic 0). 

\subsubsection*{Existence and Construction} 
\label{subsub:Exist+Construct} 

\bDef \label{UnivObj}
Let $\Mr$ be a category of rational maps from $Y$ to algebraic
groups. Then $\left(u:Y\dra\sU\right)\in\Mr$ is called
a \emph{universal object for} $\Mr$ if it has the universal mapping
property in $\Mr$: 

for all $\left(\phe:Y\dra G\right)\in\Mr$ 
there exists a unique homomorphism of algebraic groups 
$h:\sU\lra G$ such that $\phe=h\circ u$ up to translation, 
i.e.\ there is a constant $g\in G(k)$ 
such that the following diagram is commutative
\[  \xymatrix{  Y \ar@{-->}[dr]_u \ar@{-->}[rr]^{\transl_g \circ \, \phe} & & G \\  
                &  \sU \ar[ur]_h & \\ } 
\] 
where $\transl_g: x \lmt x+g$ is the translation by $g$. 
\eDef 

\bRmk Localization of a category $\Mr$ of rational maps from $Y$ to 
algebraic groups at the system of injective homomorphisms does not 
change (the equivalence class of) the universal object. 
Therefore it is often convenient to pass to the localization 
$\Insg\Mr$ ($\see$Definition \ref{H_I}). 
\eRmk 

For the category $\Mav$ of morphisms from $Y$ to abelian varieties
($\see$Definition \ref{Mav}) there exists a universal object: 
the \emph{Albanese mapping} to the \emph{Albanese variety}, 
which is denoted by $\alb:Y\lra\Alb\left(Y\right)$. 
This is a classical result (see e.g.\ \cite{La}, \cite{Ms}, \cite{S2}). 

%\vspace{\vs} 

In the following we consider categories $\Mr$ of rational maps from
$Y$ to algebraic groups satisfying the following conditions: 

\begin{tabular}{rl}
$\left( \diamondsuit \; 1 \right)$ & 
$\Mr$ contains the category $\Mav$. \\
$\left( \diamondsuit \; 2 \right)$ & 
$\left(\phe:Y\dra G\right)\in\Mr$ 
if and only if \\
 & $\forall\lambda\inn\Ld$ the induced rational map 
   $\left(\phe_{\lambda}:Y\dra\lambda_{*}G\right)\in\Mr$. \\
$\left( \diamondsuit \; 3 \right)$ & 
If $\phe:Y\dra G$ is a rational map s.t. $\rho \circ \phe$ factors through \\
 & a homomorphism $\alp: B \lra A$ of abelian varieties, then \\
 & $(\phe:Y\dra G) \in \Mr$ if and only if 
   $(\phe^{\alp}:Y\dra \alp^* G) \in \Mr$. 
\end{tabular} \\
Here $0 \ra L \ra G \overset{\rho}\ra A \ra 0$ is the canonical 
decomposition of the algebraic group $G$. 

\bThm \label{Exist univObj} 
Let $\Mr$ be a category of rational
maps from $Y$ to algebraic groups which satisfies 
$\left(\diamondsuit \; 1-3 \right)$. 
Then for $\Mr$ there exists a universal object
$\left(u:Y\dra\sU\right)\in\Mr$ if and only
if there is a formal group $\fmlG$ which is a subfunctor of $\Divf_{Y}^{0}$
such that $\Insg \Mr$ is equivalent to $\Insg \Mr_{\fmlG}$,  
where $\Mr_{\fmlG}$ is the category of rational maps which induce a 
homomorphism of formal groups to $\fmlG$ ($\see$Definition \ref{Mr_sW}). 
\eThm 

\bPf  
$\left(\Lla\right)$ Assume that $\Insg \Mr$ is equivalent
to $\Insg \Mr_{\fmlG}$, where $\fmlG$ is a formal group in $\Divf_{Y}^{0}$.
The first step is the construction of an algebraic group $\sU$ and
a rational map $u:Y\dra\sU$. In a second step the universality
of $u:Y\dra\sU$ for $\Mr_{\fmlG}$ has to be shown. 

\textbf{Step 1:} Construction of $u:Y\dra\sU$ \\
$Y$ is a regular projective variety over $k$, thus the functor $\Picf_{Y}^{0}$
is represented by an abelian variety $\Pic_{Y}^{0}$ ($\see$Corollary 
\ref{normal Pic^0}). 
Since $\fmlG \subset \Divf_Y^0$, the formal group $\fmlG$ is torsion-free. 
The natural transformation $\Divf_{Y}^{0}\lra\Picf_{Y}^{0}$
induces a 1-motive $M=\left[\fmlG\lra\Pic_{Y}^{0}\right]$. 
Let $M^{\vee}$ be the dual 1-motive of $M$. The formal group in degree $-1$ of 
$M^{\vee}$ is the Cartier-dual of the largest linear subgroup of $\Pic_Y^0$, 
and this is zero, since an abelian variety does not contain any non-trivial 
linear subgroup. Then define $\sU$ to be the algebraic group in 
degree $0$ of $M^{\vee}$, i.e.\ $\left[0\lra\sU\right]$
is the dual 1-motive of $\left[\fmlG\lra\Pic_{Y}^{0}\right]$.
The canonical decomposition 
$0\ra\sL\ra\sU\ra\sA\ra0$
is the extension of $\left(\Pic_{Y}^{0}\right)^{\vee}$ by $\fmlG^{\vee}$
induced by the homomorphism $\fmlG\lra\Pic_{Y}^{0}$ ($\see$Theorem
\ref{Struct}), where $\sL=\fmlG^{\vee}$ is the Cartier-dual of $\fmlG$
and $\sA=\left(\Pic_{Y}^{0}\right)^{\vee}$ is the dual abelian variety
of $\Pic_{Y}^{0}$, which is $\Alb\left(Y\right)$. 

As $\sL$ is a linear algebraic group,  there is a canonical splitting 
$\sL \cong \Trs \times \Vcl$ of $\sL$ into the direct product 
of a torus $\Trs$ of rank $t$ and a vectorial group $\Vcl$ of dimension $v$ 
($\see$Theorem \ref{Thm linGrp}). 
The homomorphism $\fmlG\lra\Pic_{Y}^{0}$ is uniquely determined 
by the values on a basis $\Omega$ of the finite free $\Zint$-module 
\[ \fmlG\left(k\right) = 
   \sL^{\vee} \left(k\right) = 
    \Trs^{\vee}(k) = 
   \Homabk\left(\Trs,\Gm\right)
\]  
and on a basis $\Theta$ of the finite dimensional 
$k$-vector space 
\[ \Lie\left(\fmlG\right) = 
   \Lie\left(\sL^{\vee}\right) = 
   \Lie\left(\Vcl^{\vee}\right) =
   \Homk\left(\Lie(\Vcl),k\right) = 
   \Homabk\left(\Vcl,\Ga\right) \;. 
\]  
By duality, such a choice of bases corresponds to a decomposition 
\[ \sL \overset{\sim}\lra 
   \left(\Gm\right)^{t}\times\left(\Ga\right)^{v} \;, 
\] 
and induces a decomposition 
\begin{eqnarray*}
\Ext\left(\sA,\sL\right) & \overset{\sim}\lra & 
\Ext\left(\sA,\Gm\right)^{t}\times\Ext\left(\sA,\Ga\right)^{v} \\ 
\sU & \lmt & 
\prod_{\omega\in\Omega} \omega_* \sU \times \prod_{\tha\in\Theta} \tha_* \sU \;. \\ 
\end{eqnarray*} 
Therefore the rational map $u:Y\dra\sU$ is uniquely determined by
the following rational maps to push-outs of $\sU$ 
\[ 
\begin{array}{c}
u_{\omega}:Y\dra\omega_{*}\sU\qquad\qquad\omega\in\Omega\\
u_{\vartheta}:Y\dra\vartheta_{*}\sU\qquad\qquad\vartheta\in\Theta 
\end{array} 
\] 
whenever $\Omega$ is a basis of $\fmlG\left(k\right)$
and $\Theta$ a basis of $\Lie\left(\fmlG\right)$. 
We have isomorphisms 
\begin{eqnarray*}
\Ext\left(\sA,\Gm\right)\:\simeq\:\Pic_{\sA}^{0}\left(k\right) 
& \overset{\sim}\lra & \Pic_{Y}^{0}\left(k\right)\\
P & \lmt & P_{Y}=P\times_{\sA}Y 
\end{eqnarray*} 
and 
\begin{eqnarray*}
\Ext\left(\sA,\Ga\right)\:\simeq\:\Lie\left(\Pic_{\sA}^{0}\right) 
& \overset{\sim}\lra & \Lie\left(\Pic_{Y}^{0}\right)\\
T & \lmt & T_{Y} = T\times_{\sA}Y \;. 
\end{eqnarray*}
{}From the proof of Theorem \ref{Struct} it follows that 
$\left(\omega_{*}\sU\right)_{Y}$
is just the image of 
$\omega\in\fmlG\left(k\right)\subset\Divf_{Y}^{0}\left(k\right)$
under the homomorphism $\fmlG\lra\Pic_{Y}^{0}$, which is
the divisor-class $\left[\omega\right]\in\Pic_{Y}^{0}\left(k\right)$.
Likewise from the proof of Theorem \ref{Struct} follows that 
$\left(\vartheta_{*}\sU\right)_{Y}$
is the image of 
$\vartheta\in\Lie\left(\fmlG\right)\subset\Lie\left(\Divf_{Y}^{0}\right)$
under the homomorphism $\fmlG\lra\Pic_{Y}^{0}$, which is
the class of deformation 
$\left[\vartheta\right]\in\Lie\left(\Pic_{Y}^{0}\right)$.
Then define the rational map $u:Y\dra\sU$ by the condition
that for all $\omega\in\Omega$ the section 
\[ u_{Y,\omega}:Y\dra\omega_{*}\sU_{Y} = \left[\omega\right] 
\] 
corresponds to the divisor $\omega\in\Divf_{Y}^{0}\left(k\right)$,
and for all $\vartheta\in\Theta$ the section 
\[ u_{Y,\vartheta}:Y\dra\vartheta_{*}\sU_{Y} = \left[\vartheta\right] 
\] 
corresponds to the deformation $\vartheta\in\Lie\left(\Divf_{Y}^{0}\right)$,
in the sense of Lemma \ref{div(s)}, i.e.\ \[
\begin{array}{c}
\dv_{\Gm}\left(u_{Y,\omega}\right)=\omega\qquad\qquad\forall\omega\in\Omega \;\;\\ 
\dv_{\Ga}\left(u_{Y,\vartheta}\right)=\vartheta\qquad\qquad
\forall\vartheta\in\Theta \;. 
\end{array}\] 
This determines $u$ up to translation by a constant. 
The conditions $(\diamondsuit \; 1-3)$ guarantee that $(u:Y\dra\sU) \in \Mr$. 

\textbf{Step 2:} Universality of $u:Y\dra\sU$ \\ 
Let $G$ be an algebraic group with canonical decomposition 
$0\ra L\ra G\overset{\rho}\ra A\ra0$ and $\phe:Y\dra G$ a rational map 
inducing a homomorphism of formal groups $\ld:\Ld\lra\fmlG$,
$\lambda\lmt\dv_{\splG}\left(\phe_{Y,\lambda}\right)$ for
$\lambda\inn\Ld$ ($\see$Proposition \ref{induced Trafo}). Let
$l:\sL\lra L$ be the dual homomorphism of linear groups.
The composition 
$Y\overset{\phe}\dra G\overset{\rho}\lra A$
extends to a morphism from $Y$ to an abelian variety. Translating 
$\phe$ by a constant $g\in G(k)$, if necessary, we may hence assume 
that $\rho\circ\phe$ factors through $\sA=\Alb\left(Y\right)$: 
\[  \xymatrix{  Y \ar[dr]_{\alb} \ar[rr]^{\rho\circ\phe} & & A \\  
                &  \Alb(Y) \ar[ur] & \\ } \]
We are going to show that the following diagram commutes: \[ 
\xymatrix{ && \sL \ar[dd] \ar[rr]^{l} && L \ar[dd] \ar@{=}[rr] && L \ar[dd] \\ 
        && && && \\ 
        && \sU \ar[dd] \ar[rr]^h && G_{\sA} \ar[dd] \ar[rr] && G \ar[dd]^{\rho} \\ 
        && && && \\ 
        Y \ar[rr] \ar@{-->}[uurrrrrr]^{\phe} \ar@{-->}[uurrrr]^(.45){\phe_{\sA}} 
        \ar@{-->}[uurr]^u && \sA  \ar@{=}[rr] && \sA \ar[rr] && A } \]
i.e.\ the task is to show that \[
\begin{array}{rl}
\textrm{(a)} & G_{\sA}=l_{*}\sU\\
\textrm{(b)} & \phe_{\sA}=h\circ u \qquad \textrm{mod translation} 
\end{array}\]
 where $G_{\sA}=G\times_{A}\sA$ is the fibre-product of $G$ and $\sA$
over $A$ and $\phe_{\sA}:Y\dra G_{\sA}$ is the unique
map obtained from $\left(\phe,\alb\right):Y\dra G\times\sA$
by the universal property of the fibre-product $G_{\sA}$, and where
$h$ is the homomorphism obtained by the amalgamated sum 
\[  \xymatrix{  \sL \ar[d] \ar[r] & L \ar[d] \\  
                \sU \ar[r]^-h & \sU \amalg_{\sL} L  \\ } \]
as by definition of the push-out we have $l_{*}\sU=\sU\amalg_{\sL}L$. \\
For this purpose, by additivity of extensions, it is enough to show
that for all $\lambda\inn\Ld$ it holds \[
\begin{array}{ll}
\textrm{(a')} & \lambda_{*}G_{\sA}=\ld(\lambda)_{*}\sU\\
\textrm{(b')} & \phe_{\sA,\lambda}=u_{\ld(\lambda)} \qquad \textrm{mod translation} 
\end{array}\] 
where $\ld(\lambda) = \lambda\circ l$ and 
$\ld(\lambda)_*=\left(\lambda\circ l\right)_*=\lambda_{*}l_{*}$. 
Using the isomorphism $\Pic_{\sA}^{0}\overset{\sim}\lra\Pic_{Y}^{0}$, 
this is equivalent to showing that for all $\lambda\inn\Ld$ it holds
\[
\begin{array}{ll}
\textrm{(a'')} & \lambda_{*}G_{Y}=\ld(\lambda)_{*}\sU_{Y}\\
\textrm{(b'')} & \phe_{Y,\lambda}=u_{Y,\ld(\lambda)} \qquad \textrm{mod translation} 
\;. 
\end{array} 
\] 
By construction of $u:Y\dra\sU$, we have for all $\lambda\inn\Ld$:\[
\dv_{\splG}\left(u_{Y,\ld\left(\lambda\right)}\right)=\ld(\lambda)=
\dv_{\splG}\left(\phe_{Y,\lambda}\right)\]
 and hence 
\begin{eqnarray*}
\ld(\lambda)_{*}\sU_{Y} & = & \left[\ld(\lambda)\right]\\
 & = & \left[\dv_{\splG}\left(\phe_{Y,\lambda}\right)\right]\;=\;\lambda_{*}G_{Y} \;. 
\end{eqnarray*}

As $u: Y \lra \sU$ generates $\sU$, each $h':\sU \lra G_{\sA}$ fulfilling 
$h' \circ u = \phe_{\sA}$ coincides with $h$. Hence $h$ is unique. 

$\left(\Lra\right)$ Assume that $u:Y\dra\sU$
is universal for $\Mr$. Let 
$0\ra\sL\ra\sU\ra\sA\ra0$
be the canonical decomposition of $\sU$, and let $\fmlG$ be the image
of the induced transformation $\sL^{\vee}\lra\Divf_{Y}^{0}$.
For $\lambda\inn\sL^{\vee}$  the uniqueness of the homomorphism 
$h_{\lambda}:\sU\lra\lambda_{*}\sU$
fulfilling $u_{\lambda}=h_{\lambda}\circ u$ implies that 
the rational maps $u_{\lambda}:Y\dra\lambda_{*}\sU$ are 
non-isomorphic to each other for distinct $\lambda \inn \sL^{\vee}$. 
Hence $\dv_{\splG}\left(u_{Y,\lambda}\right) \neq 
\dv_{\splG}\left(u_{Y,\lambda'}\right)$ for $\lambda \neq \lambda' \inn \sL^{\vee}$. 
Therefore $\sL^{\vee}\lra\fmlG$ is injective, hence an isomorphism. 

Let $\phe:Y\dra G$ be an object of $\Mr$ and 
$0\ra L\ra G\ra A\ra0$
be the canonical decomposition of $G$. Translating $\phe$ by
a constant $g\in G(k)$, if necessary, we may assume that 
$\phe:Y\dra G$
factorizes through a unique homomorphism $h:\sU\lra G$.
The restriction of $h$ to $\sL$ gives a homomorphism of linear groups
$l:\sL\lra L$. Then the dual homomorphism 
$\ld:\Ld\lra\fmlG$
yields a factorization of $\Ld\lra\Divf_{Y}^{0}$ through
$\fmlG$. Thus $\Mr$ is a subcategory of $\Mr_{\fmlG}$. Now the properties 
$\left(\diamondsuit \; 1-3 \right)$ guarantee that $\Insg \Mr$ contains 
the equivalence classes of all rational maps which induce a transformation 
to $\fmlG$, hence $\Insg \Mr$ is equivalent to $\Insg \Mr_{\fmlG}$.
\ePf 

\bNot  
The universal object for a category $\Mr$ of rational maps
from $Y$ to algebraic groups, if it exists, is denoted by 
$\alb_{\Mr}:Y\dra\Alb_{\Mr}\left(Y\right)$.
\\
If $\fmlG$ is a formal group in $\Divf_{Y}^{0}$, then the universal
object for $\Mr_{\fmlG}$ is also denoted by 
$\alb_{\fmlG}:Y\dra\Alb_{\fmlG}\left(Y\right)$.
\eNot  

\bRmk \label{Alb_constr} In the proof of Theorem \ref{Exist univObj}
we have seen that $\Alb_{\fmlG}\left(Y\right)$ is an extension of the
abelian variety $\Alb\left(Y\right)$ by the linear group $\fmlG^{\vee}$,
and the rational map 
$\left(\alb_{\fmlG}:Y\dra\Alb_{\fmlG}\left(Y\right)\right)\in\Mr_{\fmlG}$
is characterized by the fact that the transformation 
$\trafo_{\alb_{\fmlG}}:\Ld\lra\Divf_Y^0$ 
is the identity $\id:\fmlG\lra\fmlG$.
\\
More precisely, $\left[0\lra\Alb_{\fmlG}\left(Y\right)\right]$
is the dual 1-motive of $\left[\fmlG\lra\Pic_{Y}^{0}\right]$.
\eRmk 

\bExm \label{classical Alb} 
The universal object $\alb:Y\lra\Alb\left(Y\right)$
for $\Mav$ from Definition \ref{Mav} is the \emph{classical Albanese
mapping} and $\Alb\left(Y\right)$ the \emph{classical Albanese variety}
of a regular projective variety $Y$. \eExm 

\bExm \label{Alb_0} 
The universal object $\alb_0:Y\dra\Alb_0\left(Y\right)$ for $\Mr_0$ 
from Example \ref{Mr_0} coincides with the classical Albanese mapping 
to the classical Albanese variety by Theorem \ref{Exist univObj}, 
since $\Insg \Mr_0$ is equivalent to $\Insg \Mav$. 
\eExm

\bExm \label{Alb_0^X} 
The universal object for the category $\Mro{X}$ from Example \ref{M_0^X} 
is a quotient of the classical Albanese $\Alb \big(\Xt\big)$ of a 
projective resolution of singularities $\Xt$ for $X$, 
as $\Mro{X}$ is a subcategory of $\Mr_0$ for $Y=\Xt$. 
It is the universal object for the category of morphisms from $X$ to 
abelian varieties and 
coincides with the \emph{universal morphism for the variety $X$ 
and for the category of abelian varieties} in the sense of \cite{S2}. 
\eExm

\bExm \label{Serres genAlb} The universal object 
$\alb_{\fmlG_{D}}:Y\dra\Alb_{\fmlG_{D}}\left(Y\right)$
for $\Mr_{\fmlG_{D}}$ from Example \ref{Mr_D} is the 
\emph{generalized Albanese of Serre} (see \cite{S3}). 
\eExm 

\bExm \label{genJac} The universal object 
$\alb_{\fmlG_{\mdl}}:C\dra\Alb_{\fmlG_{\mdl}}\left(Y\right)$
for $\Mr_{\fmlG_{\mdl}}$ from Example \ref{Mr_del} is Rosenlicht's 
\emph{generalized Jacobian} $J_{\mdl}$ to the modulus $\mdl$ 
(see \cite{S}). 
\eExm 

\bExm \label{Faltings/Wuestholz} The universal object 
$\alb_{\fmlG_{W}}:Y\dra\Alb_{\fmlG_{W}}\left(Y\right)$
for $\Mr_{\fmlG_{W}}$ from Example \ref{Mr_Omega} is the 
\emph{generalized Albanese of Faltings/W\"ustholz} 
(see \cite[VI.2.]{FW}) 
\eExm

\bExm \label{univ.reg.Quot} The universal object 
$\alb_{\Divf_{\Xt/X}^{0}}:X_{\reg}\dra\Alb_{\Divf_{\Xt/X}^{0}}\big(\Xt\big)$
from Example \ref{Mr_Y/X} is the \emph{universal regular quotient}
of the Chow group of points $\CHOo{X}$ (see \cite{ESV}). 
In the following we will simply denote it by $\Alb\left(X\right)$. 
This is consistent, since in the case that $X$ is regular 
it coincides with the classical Albanese variety. 
\eExm 

\bRmk  Also the generalized Albanese of Serre ($\see$Example
\ref{Serres genAlb}) and the generalized Jacobian ($\see$Example
\ref{genJac}) 
can be interpreted as special cases of the universal regular quotient 
($\see$Example \ref{univ.reg.Quot}) 
by constructing an appropriate singular variety $X$. 
\eRmk

\subsubsection*{Functoriality} 
\label{Functoriality} 

The Question is whether a morphism of regular projective varieties induces 
a homomorphism of algebraic groups between universal objects. 

\bPrp \label{alb(sigma)} Let $\sigma: V \lra Y$ be a morphism 
of regular projective varieties. Let $\MrV$ and $\MrY$ be categories of 
rational maps from $V$ and $Y$ respectively to algebraic groups, and 
suppose there exist universal objects $\Alb_{\MrV}(V)$ and $\Alb_{\MrY}(Y)$ 
for $\MrV$ and $\MrY$ respectively. 
The universal property of $\Alb_{\MrV}(V)$ yields: \\ 
If the composition $\alb_{\MrY} \circ \,\sigma: V \dra \Alb_{\MrY}(Y)$ 
is an object of $\MrV$, then $\sigma$ induces a homomorphism of 
algebraic groups 
\[ \Alb_{\MrV}^{\MrY}(\sigma): \Alb_{\MrV}(V) \lra \Alb_{\MrY}(Y) \;. 
\] 
\ePrp

Theorem \ref{Exist univObj} allows to give a more explicit description: 

\bPrp \label{pull-back_of_1-motives} 
Let $\sigma: V \lra Y$ be a morphism of regular projective varieties. 
Let $\fmlG \subset \Divf_Y^0$ be a formal group s.t. 
$\Supp(\fmlG)$ does not contain any component of $\sigma(V)$. 
Let $\_\cut V:\Decf_{Y,V} \lra \Divf_V$ 
denote the pull-back of Cartier divisors from $Y$ to $V$ 
($\see$Definition \ref{Dec_Y,V}, Proposition \ref{_.V}). \\ 
For each formal group $\fmlGr \subset \Divf_V^0$ satisfying 
$\fmlGr \supset \fmlG\cut V$, the pull-back of relative Cartier divisors 
and of line bundles induces a transformation of 1-motives 
\[ \left[ \begin{array}{c} 
          \fmlGr \\ \downarrow \\ \Picf_V^0 \\
          \end{array} 
   \right]
   \lla 
   \left[ \begin{array}{c} 
          \fmlG \\ \downarrow \\ \Picf_Y^0 \\
          \end{array} 
   \right]
\] 
\ePrp 

Remembering the construction of the universal objects 
($\see$Remark \ref{Alb_constr}), dualization of 1-motives translates 
Proposition \ref{pull-back_of_1-motives} into the following reformulation 
of Proposition \ref{alb(sigma)}: 

\bPrp \label{alb_F(sigma)} 
Let $\sigma: V \lra Y$ be a 
morphism of regular projective varieties. Let $\fmlG \subset \Divf_Y^0$ 
be a formal group s.t. $\Supp(\fmlG)$ does not contain any component of 
$\sigma(V)$.  
Then $\sigma$ induces a homomorphism of algebraic groups 
\[ \Alb_{\fmlGr}^{\fmlG}(\sigma): \Alb_{\fmlGr}(V) \lra \Alb_{\fmlG}(Y) 
\] 
for each formal group $\fmlGr \subset \Divf_V^0$ satisfying 
$\fmlGr \supset \fmlG\cut V$. 
\ePrp

\section{Rational Maps Factoring through $\CHOo{X}$} 
\label{sec:Rat-Maps CH_0(X)_deg0} 

Throughout this section let $X$ be a projective variety over $k$ 
(an algebraically closed field of characteristic 0) 
and $\pi:Y\lra X$ a projective resolution of singularities. 
Let $U\subset Y$ be an open dense subset of $Y$ where $\pi$ is an 
isomorphism. $U$ is identified with its image in $X$, and we suppose 
$U\subset X_{\reg}$. We consider the category $\MrCH{X}$ of morphisms 
$\phe:U\lra G$ from $U$ to algebraic groups $G$ 
factoring through $\CHOo{X}$ ($\see$Definition \ref{MrCH(X)}), 
where we assume algebraic groups $G$ always to be connected, 
unless stated otherwise. 

The goal of this section is to show that the category $\MrCH{X}$ 
is equivalent to the category $\Mr_{\Divf_{Y/X}^{0}}$ of rational
maps which induce a transformation of formal groups to $\Divf_{Y/X}^{0}$, 
which is defined as the kernel in $\Divf_Y^0$ of a kind of 
push-forward $\pi_*$ of relative Cartier divisors 
(see Propositions \ref{Div_Z/C^0} and \ref{Div_Y/X^0}).

\subsection{Chow Group of Points }

In this subsection the Chow group $\CHOo{X}$
of 0-cycles of degree 0 modulo rational equivalence is presented, quite 
similar as in \cite{LW}, see also \cite{ESV}, \cite{BiS}. 

\bDef \label{Cartier-curve} A \emph{Cartier curve} in $X$, 
relative to $X\setminus U$, is a curve $C\subset X$ satisfying \\
\begin{tabular}{rl}
(a)&
$C$ is pure of dimension 1.\tabularnewline
(b)&
No component of $C$ is contained in $X\setminus U$.\tabularnewline
(c)&
If $p\in C\setminus U$, the ideal of $C$ in $\sO_{X,p}$ 
is generated by a regular sequence.\tabularnewline
\end{tabular} \eDef 

\bDef \label{R(C,X)} 
Let $C$ be a Cartier curve in $X$ relative
to $X\setminus U$, $\Cp\left(C\right)$ the set of irreducible components 
of $C$ and $\gamma_{Z}$ the generic points of $Z\in\Cp\left(C\right)$.
Let $\sO_{C,\Theta}$ be the semilocal ring on $C$ at 
$\Theta=\left(C\setminus U\right)
\cup\left\{ \gamma_{Z}\,|\, Z\in\Cp\left(C\right)\right\} $.
Define 
\[ \K\left(C,U\right)^{*} = \sO_{C,\Theta}^{*} \; . 
\] 
\eDef 

\bDef \label{div(f)_C} 
Let $C$ be a Cartier curve in $X$ relative
to $X\setminus U$ and $\nu:\Ct\lra C$ its normalization. For 
$f\in\K\left(C,U\right)^{*}$
and $p\in C$ let \[
\ord_{p}\left(f\right)=\sum_{\pt\ra p}\val_{\pt}\left(\widetilde{f}\right)\]
 where $\widetilde{f}:=\nu^{\fis}f\in\sK_{\Ct}$ and $\val_{\pt}$ is
the discrete valuation attached to the point $\pt\in\Ct$ above $p\in C$
(cf.\ \cite[Example A.3.1]{F}). \\
Define the divisor of $f$ to be 
\[ \dv\left(f\right)_{C}=\sum_{p\in C}\ord_{p}\left(f\right)\;\left[p\right] \;. 
\] 
\eDef  

\bDef \label{CH_0(X)_deg0} 
Let $\Z_{0}\left(U\right)$ be the group of 0-cycles on $U$, set \[
\fR_{0}\left(X,U\right)=\left\{ \left(C,f\right)\left|\begin{array}{c}
C\textrm{ is a Cartier curve in }X\textrm{ relative to } X\setminus U\\
\textrm{and }f\in\K\left(C,U\right)^{*}\end{array}\right.\right\} \]
 and let $\R_{0}\left(X,U\right)$ be the subgroup of $\Z_{0}\left(U\right)$
generated by the elements $\dv\left(f\right)_{C}$ with 
$\left(C,f\right)\in\fR_{0}\left(X,U\right)$.
Then define \[
\CHO{X}=\Z_{0}(U)/\R_{0}(X,U) \;. \] 
Let $\CHOo{X}$ be the subgroup of $\CHO{X}$
of cycles $\zeta$ with $\deg\zeta|_{W}=0$ for all irreducible components
$W\in\Cp\left(U\right)$ of $U$. \eDef  

\bRmk  
The definition of $\CHO{X}$ and $\CHOo{X}$
is independent of the choice of the dense open subscheme $U\subset X_{\reg}$
(see \cite[Corollary 1.4]{ESV}). 
\eRmk  

\bRmk  
Note that by our terminology a curve is always reduced, in
particular a Cartier curve. In the literature, e.g.\ \cite{ESV}, \cite{LW},
a slightly different definition of \emph{Cartier curve} seems to be
common, which allows non-reduced Cartier curves. Actually this does
not change the groups $\CHO{X}$ and $\CHOo{X}$, 
see \cite[Lemma 1.3]{ESV} for more explanation. 
\eRmk

\subsection{Local Symbols} 

The description of rational maps factoring through $\CHOo{X}$ requires 
the notion of a \emph{local symbol} as in \cite[Chapter III, \S 1]{S}. 

Let $C$ be a smooth projective curve over $k$. 
The composition law of an unspecified algebraic group $G$ 
is written additively in this subsection. 

\bDef  
For an effective divisor $\mdl=\sum n_{p}\, p$ on $C$, 
a subset $S\subset C$ and rational functions $f,g\in\sK_{C}$ define 
\begin{eqnarray*} 
f\equiv g\mod\mdl\textrm{ at }S & \quad:\Llra\quad & 
   \val_{p}\left(f-g\right)\geq n_{p}
   \qquad\forall p\in S\cap\Supp\left(\mdl\right) \; ,\\ 
f\equiv g\mod\mdl & \quad:\Llra\quad & 
   f\equiv g\mod\mdl\textrm{ at }C \; . 
\end{eqnarray*}
where $\val_{p}$ is the valuation attached to the point $p\in C$. 
\eDef 

Let $\psi:C\dra G$ be a rational map from $C$ to an algebraic
group $G$ which is regular away from a finite subset $S$. The morphism
$\psi:C\setminus S\lra G$ extends to a homomorphism from
the group of 0-cycles $\Z_{0}\left(C\setminus S\right)$ to $G$ by
setting $\psi\left(\sum l_i\, c_i\right):=\sum l_i\,\psi\left(c_i\right)$
for $c_i\in C\setminus S$, $l_i\in\Zint$. 

\bDef \label{modulus} 
An effective divisor $\mdl$ on $C$ is said
to be a \emph{modulus} for $\psi$ if $\psi\left(\dv\left(f\right)\right)=0$
for all $f\in\sK_{C}$ with $f\equiv1\mod\mdl$. \eDef 

\bThm \label{modulus_exist} 
Let $\psi:C\dra G$ be a
rational map from $C$ to an algebraic group $G$ and $S$ the finite
subset of $C$ where $\psi$ is not regular. Then $\psi$ has a modulus
supported on $S$. \eThm 

This theorem is proven in \cite[Chapter III, \S 2]{S}, using the
following concept: 

\bDef \label{locSymbol} 
Let $\mdl$ be an effective divisor supported
on $S\subset C$ and $\psi:C\dra G$ a rational function
from $C$ to an algebraic group $G$, regular away from $S$. A \emph{local
symbol} associated to $\psi$ and $\mdl$ is a function \[
\left(\psi,\_\right)_{\_}:\sK_{C}^{*}\times C\lra G\]
 which assigns to $f\in\sK_{C}^{*}$ and $p\in C$ an element 
$\left(\psi,f\right)_{p}\in G$,
satisfying the following conditions: 

\begin{tabular}{rl}
(a)&
$\left(\psi,fg\right)_{p}=\left(\psi,f\right)_{p}+\left(\psi,g\right)_{p}$ \, , 
\tabularnewline
(b)&
$\left(\psi,f\right)_{c}=\val_{c}(f)\;\psi(c)\qquad$ 
if
$c\in C\setminus S$ \ , \tabularnewline
(c)&
$\left(\psi,f\right)_{s}=0\qquad\qquad$if $s\in S$ and
$f\equiv1\mod\mdl$ at $s$ \ , \tabularnewline
(d)&
$\sum_{p\in C}\left(\psi,f\right)_{p}=0$ \, . \tabularnewline
\end{tabular} \eDef 

\bPrp \label{locSymb_unique} The rational map $\psi$ has a modulus
$\mdl$ if and only if there exists a local symbol associated to $\psi$
and $\mdl$, and this symbol is then unique. \ePrp 

\bPf  \cite[Chapter III, No.~1, Proposition 1]{S}. 
\ePf  

\vspace{\vs}  

Theorem \ref{modulus_exist} in combination with Proposition 
\ref{locSymb_unique}
states for each rational map $\psi:C\dra G$ the existence
of a modulus $\mdl$ for $\psi$ and of a unique local symbol 
$\left(\psi,\_\right)_{\_}$
associated to $\psi$ and $\mdl$. 

{}From the definitions it is clear that if $\mdl$ is a modulus for
$\psi$ then $\mdll$ is also for all $\mdll\geq\mdl$. Likewise a
local symbol $\left(\psi,\_\right)_{\_}$ associated to $\psi$ and
$\mdl$ is also associated to $\psi$ and $\mdll$ for all $\mdll\geq\mdl$. 

Suppose we are given two moduli $\mdl$ and $\mdl'$ for $\psi$,
and hence two local symbols $\left(\psi,\_\right)_{\_}$ and 
$\left(\psi,\_\right)'_{\_}$
associated to $\mdl$ and $\mdl'$ respectively. Then both local symbols
are also associated to $\mdll:=\mdl+\mdl'$. The uniqueness of the
local symbol associated to $\psi$ and $\mdll$ implies that 
$\left(\psi,\_\right)_{\_}$
and $\left(\psi,\_\right)'_{\_}$ coincide. It is therefore morally
justified to speak about \emph{the local symbol associated to} $\psi$
(without mentioning a modulus), 
cf.\ \cite[Chapter III, No.~1, Remark of Proposition 1]{S}. 

\bCor \label{locSymb_fml} For each rational map $\psi:C\dra G$
from $C$ to an algebraic group $G$ there exists a unique associated
local symbol $\left(\psi,\_\right)_{\_}:\sK_{C}^{*}\times C\lra G$.
If $\mdl$ is a modulus for $\psi$ supported on $S$, then this local
symbol is given by \begin{eqnarray*}
\left(\psi,f\right)_{c} & = & \val_{c}(f)\;\;\psi(c)
\qquad\qquad\qquad\forall c\in C\setminus S\\
\left(\psi,f\right)_{s} & = & 
-\sum_{c\notin S}\val_{c}(f_{s})\;\;\psi(c)
\qquad\qquad\forall s\in S 
\end{eqnarray*}
 where $f_{s}\in\sK_{C}^{*}$ is a rational function with $f_{s}\equiv1\mod\mdl$
at $z$ for all $z\in S\setminus s$ and $f/f_{s}\equiv1\mod\mdl$
at $s$. \eCor  

The above formula is shown in \cite[Chapter III, No.~1]{S}, in the
proof of Proposition 1. 

\bExm \label{locSymb_Gm} In the case that $G$ is the multiplicative
group $\Gm$, a rational map $\psi:C\dra\Gm$ can be identified
with a rational function in $\sK_{C}$, and $S$ is the set of zeros
and poles of $\psi$, i.e.\ $S=\Supp\left(\dv\left(\psi\right)\right)$.
Then the local symbol associated to $\psi$ is given by \[
\left(\psi,f\right)_{p}=\left(-1\right)^{mn}\frac{\psi^{m}}{f^{n}}(p)\qquad\;
\textrm{with }\, m=\val_{p}\left(f\right),\,n=\val_{p}\left(\psi\right) \;. 
\] 
(See \cite[Chapter III, No.~4, Proposition 6]{S}.) \eExm 

\bExm \label{locSymb_Ga} In the case that $G$ is the additive group
$\Ga$, a rational map $\psi:C\dra\Ga$ can be identified
with a rational function in $\sK_{C}$, and $S$ is the set of poles
of $\psi$. Then the local symbol associated to $\psi$ is given by
\[
\left(\psi,f\right)_{p}=\Res_{p}\left(\psi\;\der f/f\right) \;. 
\] 
(See \cite[Chapter III, No.~3, Proposition 5]{S}.) \eExm 

\bPrp \label{locSymb_additiv} Let $\phe,\psi:C\dra G$
be two rational maps from $C$ to an algebraic group $G$, with associated
local symbols $\left(\phe,\_\right)_{\_}$ and $\left(\psi,\_\right)_{\_}$.
Then the local symbol $\left(\phe+\psi,\_\right)_{\_}$ associated
to the rational map $\phe+\psi:C\dra G$, 
$c\lmt\phe\left(c\right)+\psi\left(c\right)$
is given by \[
\left(\phe+\psi,f\right)_{p}=
\left(\phe,f\right)_{p}+\left(\psi,f\right)_{p} \;. 
\] 
\ePrp  

\bPf  Let $\mdl_{\phe}$ be a modulus for $\phe$ and $\mdl_{\psi}$
one for $\psi$. Then both maps $\phe$, $\psi$ and the map $\phe+\psi$
have $\mdl_{\phe+\psi}:=\mdl_{\phe}+\mdl_{\psi}$ as a modulus
and both local symbols $\left(\phe,\_\right)_{\_}$ and 
$\left(\psi,\_\right)_{\_}$
are associated to $\mdl_{\phe+\psi}$. Now the formula in Corollary
\ref{locSymb_fml} and the distributive law imply the assertion. \ePf  

\bLem \label{locSymb_triv} 
Let $\psi:C\dra G$ be a rational map from $C$ to an algebraic group $G$ 
which is an $L$-bundle over an abelian variety $A$, 
i.e.\ $G\in\Ext\left(A,L\right)$,
where $L$ is a linear group. Let $p\in C$ be a point, $U\ni p$
a neighbourhood and $\Phi:U\times L\overset{\sim}\lra G_{U}$,
$\left(u,l\right)\lmt\phi\left(u\right)+l$ a local trivialization
of the induced $L$-bundle $G_{C}=G\times_{A}C$ over $C$, i.e.\ 
$\phi:U\lra G_{C}$
a local section. Moreover let $\left[\psi\right]_{\Phi}:C\dra L$,
$c\lmt\psi\left(c\right)-\phi\left(c\right)$ be the rational
map $\psi$ considered in the local trivialization $\Phi$. Then for
each rational function $f\in\sO_{C,p}^{*}$ 
it holds 
\[ \left(\psi,f\right)_{p}=\left(\left[\psi\right]_{\Phi},f\right)_{p} \;. 
\] 
\eLem  

\bPf  Proposition \ref{locSymb_additiv} yields 
\[ 
\left(\left[\psi\right]_{\Phi},f\right)_{p} = \left(\psi-\phi,f\right)_{p} 
 = \left(\psi,f\right)_{p}-\left(\phi,f\right)_{p} \;. 
\] 
$\phi$ is regular at $p$, therefore we have 
$\left(\phi,f\right)_{p} = \val_{p}\left(f\right)\cdot\phi\left(p\right)$. 
Since $f$ is a unit at $p$, %$f\in\sO_{C,p}^{*}$ 
it holds $\val_{p}\left(f\right)=0$. 
Thus $\left(\phi,f\right)_{p} = 0$ and hence 
$\left(\left[\psi\right]_{\Phi},f\right)_{p} = \left(\psi,f\right)_{p}$. 
\ePf

\subsection{Formal Infinitesimal Divisors} 

For a $k$-scheme $Y$ 
the functor of relative Cartier divisors $\Divf_{Y}$ admits a pull-back,
but not a push-forward. %of relative Cartier divisors. 
Supposed $Y$ is a normal scheme, the group of Cartier divisors $\Div(Y)$ 
on $Y$ can be identified with the group of locally principal Weil divisors, 
and there is a push-forward of Weil divisors. 

We are looking for a concept of infinitesimal divisors $\LDiv(Y)$
which admits a push-forward and a transformation 
$\Lie\left(\Divf_{Y}\right) \lra \LDiv(Y)$.  
In this subsection we consider the case that $Y$ is a curve $\Curv$.

\subsubsection*{Functor of Formal Infinitesimal Divisors} 

Let $\Curv$ be a curve over $\fld$. 

\bDef \label{LDiv} 
Define the $\fld$-vector space of 
\emph{formal infinitesimal divisors on $\Curv$} by 
\[ \LDiv(\Curv) = \bigoplus_{\pnt\in\Curv(\fld)} 
   \Hom_{\fld}^{\cont} \left( \mc_{\Curv,\pnt},\fld \right)
\] 
where $\Hom_{\fld}^{\cont}$ denotes the set of continuous $\fld$-linear maps. 
$\mc_{\Curv,\pnt}$ carries the $\mc_{\Curv,\pnt}$-adic topology, 
while $\fld$ is endowed with the discrete topology. 
\eDef 

\bPrp \label{phi_*} 
Let $\pi:\Curv\lra\Crvv$ be a finite morphism of curves over $\fld$. 
Then $\pi$ induces a \emph{push-forward} of formal infinitesimal divisors 
\[ \pi_{*}:\LDiv(\Curv)\lra\LDiv(\Crvv) 
\] 
induced by the homomorphisms 
\[ \Hom_{\fld}^{\cont} \left( \mc_{\Curv,\pnt},\fld \right)  \lra  
  \Hom_{\fld}^{\cont} \left( \mc_{\Crvv,\pi(\pnt)},\fld \right), \quad 
  h  \lmt  h \circ \widehat{\pi^{\fis}} 
\] 
where $\pnt\in\Curv(\fld)$ and 
$\widehat{\pi^{\fis}}: \Oc_{\Crvv,\pi(\pnt)} \lra \Oc_{\Curv,\pnt}$ 
is the homomorphism of completed structure sheaves associated to $\pi$. 
\ePrp 

\bPrp \label{fml} 
Let $\Curv$ be a normal curve over $\fld$. 
Then there is an isomorphism of $\fld$-vector spaces  
\[ \fml: \Lie\left(\Divf_{\Curv}\right)  \lra  \LDiv(\Curv) \;. 
\] 
\ePrp 

\bPf 
We construct the isomorphism $\fml$ via factorization, 
i.e.\ give isomorphisms 
\[ \Gam\left(\left.\sK_{\Curv}\right/\sO_{\Curv}\right) \overset{\sim}\lra 
   \bigoplus_{\pnt\in\Curv(\fld)} \left.\sK_{\Curv,\pnt}\right/\sO_{\Curv,\pnt} 
   \overset{\sim}\lra 
   \bigoplus_{\pnt\in\Curv(\fld)} \Hom_{\fld}^{\cont}\left(\mc_{\Curv,\pnt},\fld\right) \;. 
\] 
The first of these two maps is given by the natural $\fld$-linear map 
\[ \Gamma\left(\left.\sK_{\Curv}\right/\sO_{\Curv}\right)  \lra 
   \bigoplus_{\pnt\in \Curv(\fld)} \left.\sK_{\Curv,\pnt}\right/\sO_{\Curv,\pnt}, \quad
   \delta  \lmt  \sum_{\pnt\in \Curv(\fld)} [\delta]_{\pnt} 
\] 
which is an isomorphism by the Approximation Lemma 
(see \cite[Part One, Chapter I, \S 3]{Se}). 
As $\Curv$ is normal, each local ring is regular. 
Since 
$\left.\sK_{\Curv,\pnt}\right/\sO_{\Curv,\pnt} = 
   \bigcup_{\nu>0} \left. t_{\pnt}^{-\nu}\,\Oc_{\Curv,\pnt} \right/ \Oc_{\Curv,\pnt} 
$ 
for a local parameter $t_{\pnt}$ of the maximal ideal 
$\fm_{\pnt}\subset\sO_{\Curv,\pnt}$, 
Lemma \ref{l cont} below yields a canonical isomorphism 
of $\fld$-vector spaces 
%\[ \bigcup_{\nu>0}\left.t_{\pnt}^{-\nu}\,\Oc_{\Curv,\pnt} \right/ 
%    \Oc_{\Curv,\pnt} 
\[  \left.\sK_{\Curv,\pnt}\right/\sO_{\Curv,\pnt} 
    \overset{\sim}\lra  \Hom_{\fld}^{\cont}\left(\mc_{\Curv,\pnt},\fld\right), \quad 
    [f] \lmt \Res_{\pnt}\left(f\cdot\der\_\right) \;. 
\] 
Then the isomorphism $\fml$ is obtained by composition. 
\ePf 

\bLem \label{l cont} 
Let $(\sA,\fm)$ be a complete local $\fld$-algebra, 
endowed with the $\fm$-adic topology, 
while $\fld$ carries the discrete topology. 
Let $l\in\Hom_{\fld}\left(\fm,\fld\right)$ be a $\fld$-linear map. 
Then the following conditions are equivalent: 

\begin{tabular}{rl}
(i)  & $l$ is continuous, \\
(ii) & $\ker\left(l\right)$ is open, \\
(iii)& $\ker\left(l\right)\supset\fm^{\nu}$ for some $\nu>0$, \\
(iv) & $l\in\Hom_{\fld}\left(\fm/\fm^{\nu},\fld\right)$ for some $\nu>0$. \\ 
\end{tabular} \\ 
If furthermore $\sA$ is a discrete valuation ring, this is equivalent to 

\begin{tabular}{rl}
(v)&$l=\Res\left(f\cdot\der\_\right):g\lmt\Res\left(f\cdot\der g\right)$
        for some $f\in t^{-\nu}\sA/\sA,\;\nu\geq0$ \\
\phantom{(iii)} & where $t$ is a local parameter of $\fm$, 
   $\sK = \Q\left(\sA\right)$ the quotient field, \\ %of $\sA$ \\ 
 & $\Res:\Omega_{\sK/\fld}\lra\fld$ the residue and 
   $\der:\sA\lra\Omega_{\sA/\fld}$ the universal \\ 
 & derivation. \\ 
\end{tabular}
\eLem 

\bPf  
(i)$\Llra$(ii)$\Llra$(iii) is folklore 
of rings with $\fm$-adic topology. \\
(iii)$\Llra$(iv) 
$\Hom_{\fld}\left(\fm/\fm^{\nu},\fld\right) = 
 \ker\big(\Hom_{\fld}\left(\fm,\fld\right) \lra 
 \Hom_{\fld}\left(\fm^{\nu},\fld\right)\big)$. \\
(iv)$\Llra$(v) If $\sA$ is regular and $t \in \fm$ a local parameter, 
then we may identify $\sA \cong \fld[[t]]$ and $\sK \cong \fld((t))$. 
The residue over $\fld$ is defined as 
\[ 
\Res:\;\Omega_{\sK/\fld}  \lra  \fld, \quad 
\sum_{\nu\gg -\infty} a_{\nu} t^{\nu} \;\der t  \lmt  a_{-1}
\]
and the definition is independent of the choice of local parameter, 
(see \cite[Chapter II, No.~7, Proposition 5]{S}). 

$\der:\sA\lra\Omega_{\sA/\fld}$ and $\Res:\Omega_{\sK/\fld}\lra\fld$ are both 
$\fld$-linear maps. Since $\Res\left(\omega\right)=0$ for all 
$\omega\in\Omega_{\sA/\fld}$, the expression $\Res\left(f\;\der g\right)$ 
is well defined for $g\in\fm/\fm^{\nu+1}$ and $f\in t^{-\nu}\sA/\sA$. 

The pairing $\quad t^{-\nu}\sA/\sA\times\fm/\fm^{\nu+1}\lra
\fld,\quad\left(f,g\right)\lmt\Res\left(f\;\der g\right)\quad$
is a perfect pairing, hence 
$ t^{-\nu}\sA/\sA\overset{\sim}\lra
\Hom_{\fld}\left(\fm/\fm^{\nu+1},\fld\right),\quad f\lmt
\Res\left(f\;\der\_\right)$
is an isomorphism. 
\ePf  

\bLem \label{IDiv_Y/X^0 repr} 
Let $\Crvv$ be a projective curve over a field $\fld$, 
and let $\pi:\Curv\lra\Crvv$ be its normalization. 
Then the kernel of the composition $\pi_* \circ \fml$ 
\[ \ker \bigg( \Lie\left( \Divf_{\Curv} \right) 
               \overset{\fml}\lra \LDiv(\Curv) 
               \overset{\pi_*}\lra \LDiv(\Crvv) 
        \bigg) 
\] 
is a finite dimensional $\fld$-vector space. 
More precisely, if $S$ denotes the inverse image in $\Curv$ of the 
singular locus of $\Crvv$, for each $\pnt \in S$ there is an integer 
$n_{\pnt} \geq 0$ such that 
$\dsum_{\pnt\ra\pntt} \mc_{\Curv,\pnt}^{n_{\pnt}+1} \subset \mc_{\Crvv,\pntt}$. 
Then 
\[ \ker\left(\pi_* \circ \fml\right) \; \subset \; 
   \Gamma\left(
       \left.\sO_{\Curv}\left(\sum_{\pnt\in S} n_{\pnt}\,\pnt\right)\right/\sO_{\Curv}
             \right) 
\]
\eLem 

\bPf 
Since the normalization is birational, the set $S$ of
$\fld$-rational points $\pnt\in\Curv$ such that 
$ \left(\sO_{\Curv,\pnt},\fm_{\Curv,\pnt}\right) \neq 
  \left(\sO_{\Crvv,\pi(\pnt)},\fm_{\Crvv,\pi(\pnt)}\right) $
is finite. 
As $\pi_*\sO_{\Curv} / \sO_{\Crvv} = 
  \prod_{\pntt\in\pi(S)} \sO_{\Curv,\pntt} / \sO_{\Crvv,\pntt}$ 
is a coherent sheaf, 
$\sO_{\Curv,\pntt} / \sO_{\Crvv,\pntt}$ 
is finite dimensional for each $\pntt\in\pi(S)$, 
hence compatible with completion. 
Thus 
$\left.\left(\dsum_{\pnt\ra\pntt} \mc_{\Curv,\pnt}\right) \right/ \mc_{\Crvv,\pntt} 
\subset \Oc_{\Curv,\pntt} / \Oc_{\Crvv,\pi(\pntt)} 
= \sO_{\Curv,\pntt} / \sO_{\Crvv,\pi(\pntt)}$ 
is also finite dimensional. 
We obtain 
\begin{eqnarray*}
 &   & \ker\Big(\LDiv(\Curv) \lra \LDiv(\Crvv) \Big) \\
 & = & \ker \Bigg( 
\dsum_{\pnt\in\Curv(\fld)} \Hom_{\fld}^{\cont}\Big(\mc_{\Curv,\pnt},\fld\Big) \lra 
\dsum_{\pntt\in\Crvv(\fld)} \Hom_{\fld}^{\cont}\Big(\mc_{\Crvv,\pntt},\fld\Big) 
            \Bigg) \\
 & = & \dsum_{\pntt\in\Crvv(\fld)} \ker \Bigg( 
\dsum_{\pnt\ra\pntt} \Hom_{\fld}^{\cont}\Big(\mc_{\Curv,\pnt},\fld\Big) \lra 
\Hom_{\fld}^{\cont}\Big(\mc_{\Crvv,\pntt},\fld\Big)
                                \Bigg) \\
 & = & \dsum_{\pntt\in\Crvv(\fld)} \ker \Bigg( 
\Hom_{\fld}^{\cont}\bigg(\dsum_{\pnt\ra\pntt} \mc_{\Curv,\pnt},\fld\bigg) \lra 
\Hom_{\fld}^{\cont} \Big( \mc_{\Crvv,\pntt} , \fld \Big)
                                   \Bigg) \\
 & = & \dsum_{\pntt\in\pi(S)} \Hom_{\fld}\Bigg(
\bigg.\bigg(\dsum_{\pnt\ra\pntt} \mc_{\Curv,\pnt}\bigg) \bigg/ 
\mc_{\Crvv,\pntt} \;,\; \fld \Bigg) 
\end{eqnarray*} 
is finite dimensional. 
Since \;$\fml: \Lie\left(\Divf_{\Curv}\right)  \lra  \LDiv(\Curv)$ is injective 
by Proposition \ref{fml}, 
it follows that \;$\ker\left(\pi_* \circ \fml\right)$\; is finite dimensional. 

The finiteness of the dimension of 
$\left.\left(\dsum_{\pnt\ra\pntt} \mc_{\Curv,\pnt}\right) \right/ \mc_{\Crvv,\pntt}$
implies that for each $\pnt \in S$ there is an integer $n_{\pnt} \geq 0$ 
such that 
$\dsum_{\pnt\ra\pntt} \mc_{\Curv,\pnt}^{n_{\pnt}+1} \subset \mc_{\Crvv,\pntt}$. 
Then 
\begin{eqnarray*} 
\ker\big(\LDiv(\Curv) \lra \LDiv(\Crvv) \big) & \subset & 
\dsum_{\pnt\in S} \Hom_{\fld}\left(\mc_{\Curv,\pnt}
        \left/\left(\mc_{\Curv,\pnt}\right)^{n_{\pnt}+1} \right.,\fld\right) \;. 
\end{eqnarray*} 
If $t_{\pnt}$ is a local parameter of $\mc_{\Curv,\pnt}$, 
Lemma \ref{l cont} (iv)$\Llra$(v) yields 
\begin{eqnarray*}
   \Hom_{\fld}\left(\mc_{\Curv,\pnt}
        \left/\left(\mc_{\Curv,\pnt}\right)^{n_{\pnt}+1} \right.,\fld\right) 
   & \cong & \dsum_{\pnt\in S}t_{\pnt}^{-n_{\pnt}}\Oc_{\Curv,\pnt}/\Oc_{\Curv,\pnt} 
\end{eqnarray*} 
Then 
\begin{eqnarray*} 
\ker\left(\pi_* \circ \fml\right) & \subset & \fml^{-1}
   \left(\dsum_{\pnt\in S} t_{\pnt}^{-n_{\pnt}}\Oc_{\Curv,\pnt}/\Oc_{\Curv,\pnt}\right) \\
 & = & \Gamma\left(
       \left.\sO_{\Curv}\left(\sum_{\pnt\in S} n_{\pnt}\,\pnt\right)\right/\sO_{\Curv}
             \right) \; . 
\end{eqnarray*} 
\ePf

\subsection{The Functor $\Divf_{Y/X}^{0}$} 
\label{sub: Div_Y/X^0}

The idea about $\Divf_{Y/X}^{0}$ is to define a functor which admits
a natural transformation to the Picard functor $\Picf_{Y}^{0}$ and
measures the difference between the schemes $Y$ and $X$, 
where $\pi:Y\lra X$ is a projective resolution of singularities. 
Roughly speaking, $\Divf_{Y/X}^{0}$ is a subfunctor 
of $\Divf_{Y}^{0}$ which lies in the kernel of some kind of 
push-forward $\pi_{*}$. 

\bDef \label{weil} 
For a $k$-scheme $Y$ denote by $\WDiv(Y)$ the abelian group 
of Weil divisors on $Y$. 
Write $\weil:\Div(Y)\lra\WDiv(Y)$ for the homomorphism which 
maps a Cartier divisor to its associated Weil divisor, 
as defined in \cite[2.1]{F}. 
\eDef 

\bPrp \label{Div_Z/C^0} 
Let $\Crvv$ be a projective curve over $k$, 
and let $\pi:\Curv\lra\Crvv$ be its normalization. 
Then there is a subfunctor $\Divf_{\Curv/\Crvv}^{0}$ of $\Divf_{\Curv}^0$, 
represented by a formal group, characterized by the following conditions: 
\begin{eqnarray*}
  \Divf_{\Curv/\Crvv}^{0}(k) & = & 
     \ker\left(\Divf_{\Curv}^{0}(k)\overset{\weil}\lra
     \WDiv(\Curv)\overset{\pi_{*}}\lra\WDiv(\Crvv)\right) \\ 
  \Lie\left(\Divf_{\Curv/\Crvv}^{0}\right) & = & 
     \ker\left(\Lie\left(\Divf_{\Curv}^{0}\right)\overset{\fml}\lra
     \LDiv(\Curv)\overset{\pi_{*}}\lra\LDiv(\Crvv)\right)  \;. 
\end{eqnarray*}
\ePrp 

\bPf 
A formal group in characteristic 0 is determined by its $k$-valued 
points and its Lie-algebra ($\see$Corollary \ref{fmlG determin}). 
Then the conditions on $\Divf_{\Curv/\Crvv}^{0}$ 
determine uniquely  a subfunctor of $\Divf_{\Curv}$ 
(cf.\ Proposition \ref{fmlG->Div}). 
Thus it suffices to show that $\Divf_{\Curv/\Crvv}^0(k)$ is a free abelian group 
of finite rank and $\Lie\big(\Divf_{\Curv/\Crvv}^{0}\big)$ is 
a $k$-vector space of finite dimension. 
The latter assertion was proven in Lemma \ref{IDiv_Y/X^0 repr}. 
For the first assertion note that the normalization $\pi:\Curv\lra\Crvv$ 
is an isomorphism on the regular locus of $\Crvv$. As $\Curv$ is normal, 
$\weil:\Divf_{\Curv}(k)\lra\WDiv(\Curv)$ is an isomorphism. 
Then $\Divf_{\Curv/\Crvv}^0(k)$ is contained in the free abelian group generated 
by the preimages of the singular points of $\Crvv$, of which there exist 
only finitely many. $\Divf_{\Curv/\Crvv}^0(k)$ being a subgroup of a finitely 
generated free abelian group is also free abelian of finite rank. 
\ePf 

\bPrp \label{Div_Y/X^0} 
Let $X$ be a projective variety over $k$, 
and let $\pi: Y \lra X$ be a projective resolution of singularities. 
Let $\fmlF: \Algk \lra \Ab$ be the functor 
\[  \fmlF \; = \; \bigcap_{C} \left( \_\cut\Ct \right)^{-1} \Divf_{\Ct/C}^0 
\] 
where $C$ ranges over all Cartier curves in $X$ relative to the singular 
locus $X_{\sing}$ ($\see$Definition \ref{Cartier-curve}), 
$\Ct$ is the normalization of $C$ 
and $\_\cut \Ct:\Decf_{Y,\Ct}\lra\Divf_{\Ct}$ 
the pull-back of relative Cartier divisors from $Y$ to $\Ct$ 
($\see$Definition \ref{Dec_Y,V}, Proposition \ref{_.V}). 

Then there is a subfunctor $\Divf_{Y/X}^{0}$ of $\Divf_{Y}^{0}$, 
represented by a formal group, characterized by the conditions 
$\Divf_{Y/X}^{0}(k) = \fmlF(k)$ and 
$\Lie\big(\Divf_{Y/X}^{0}\big) = \Lie\left(\fmlF\right)$. 
\ePrp 

\bRmk \label{val} 
Let $\delta\in\Lie\lp \Divf_Y^0\rp = 
\Gamma \left(\left.\sK_{Y}\right/\sO_{Y}\right)$ be a 
deformation of the zero divisor in $Y$. Then $\delta$ determines an 
effective divisor by the poles of its local sections. 
Hence for each generic point $\eta$ of height 1 in $Y$, with 
associated discrete valuation $\val_{\eta}$, 
the expression $\val_{\eta}(\delta)$ is well defined and 
$\val_{\eta}(\delta)\leq 0$.  %if $\del\neq 0$ at $\eta$. 
Thus we obtain a homomorphism 
$\val_{\eta}:\Lie\lp \Divf_Y^0\rp \lra \Zint$.  %\cup \{\infty\}$. 
\eRmk 

\bPf[Proof of Prop.~\ref {Div_Y/X^0}] 
As in the proof of Proposition \ref{Div_Z/C^0} it suffices to 
show that $\Divf_{Y/X}^{0}\left(k\right)$ is a free abelian group of finite 
rank and $\Lie\big(\Divf_{Y/X}^{0}\big)$ is a $k$-vector space of finite 
dimension. 

Let $D \in \Divf_Y\left(k\right)$ be a non-trivial divisor on $Y$ 
whose support is not contained in the inverse image 
$\SgY=\Sing\tms_X Y$ of the singular locus $\Sing=X_{\sing}$ of $X$. 
Then $\pi\big(\Supp(D)\big)$ on $X$ is not contained in $\Sing$. 
Let $\sL$ be a very ample line bundle on $X$, 
consider the space $\left|\sL\right|^{d-1}$, where $d = \dim X$, 
of complete intersection curves $C = H_1 \cap \ldots \cap H_{d-1}$ 
with $H_i \in \left|\sL\right| = \Prj\big(\H^0(X,\sL)\big)$ 
for $i = 1,\ldots,d-1$. 
For Cartier curves $C$ in $\left|\sL\right|^{d-1}$ 
the following properties are open and dense: 

\begin{tabular}{rl} 
(a) & $C_Y = C \tms_X Y$ is regular. \\ 
(b) & $C$ intersects $\pi\big(\Supp(D)\big) \cap X_{\reg}$ properly. \\ 
(c) & $D \cut C_Y$ is a non-trivial divisor on $C_Y \cap X_{\reg}$. 
\end{tabular} \\ 
(a) is a consequence of the Bertini theorems, 
(b) is due to the fact that $\sL$ is very ample and 
(c) follows from (b) and the fact that $\Supp(D)$ 
is locally a prime divisor at almost every $q\in Y$. 
Therefore there exists a Cartier curve $C$ in $X$ 
satisfying the conditions (a)-(c). 
Then the normalization $\nu:\Ct\lra C$ coincides with $\pi|_{C_Y}$ 
and hence is an isomorphism on $C_Y \cap X_{\reg}$. 
Thus $\nu_*\big(D\cut\Ct\big) \neq 0$. 
This implies $D \notin \Divf_{Y/X}^{0}\left(k\right)$. 
Hence $\Divf_{Y/X}^{0}\left(k\right)$ is a subgroup of the free abelian group 
generated by the irreducible components of $\SgY$ of codimension 1. 
As $\SgY$ has only finitely many components, this group has finite rank. 
So $\Divf_{Y/X}^{0}\left(k\right)$ is a subgroup of a free abelian group of 
finite rank, hence is also free abelian of finite rank. 

Now let $\del \in \Lie\left(\Divf_Y\right)$ be a deformation of the trivial 
divisor on $Y$. The same argument as above shows that if 
$\del \in \Lie\big(\Divf_{Y/X}^{0}\big)$, then $\Supp(\del) \subset \SgY$. 
If $C$ is a Cartier curve in $X$ relative to $X_{\sing}$, 
we denote by $C^Y$ the proper transform of $C$, 
i.e.\ the closure of $\pi^{-1}(C \isec X_{\reg})$ in $Y$. 
As $\pi|_{C^Y}:C^Y \lra C$ is a birational morphism, 
the normalization $\nu: \Ct \lra C$ 
factors through a morphism $\mu: \Ct \lra C^Y$. 
Given $\eta\in\SgY^{\htl}$,  %\cap \Supp(\del)$, 
where $\SgY^{\htl}$ denotes the set of generic points of $\SgY$ 
of height 1 in $Y$, 
let $Z$ be a curve in $Y$ which 
intersects the prime divisor $\Es_{\eta}$ associated to $\eta$ 
properly in a point $p$.  
As $Y$ is regular, $\sO_{Y,\eta}$ is a discrete valuation ring. 
Let $\val_{\eta}$ be the valuation at $\eta$, 
and let $\val_q$ be the valuation 
attached to a point $q \in \Ct$ above $p \in C_Y$. 
Since $\Lie\big(\Divf_{\Ct/C}^{0}\big)$ is finite dimensional, 
there exists a number $n_q \in \Nat$ such that 
$\val_q (\gam) \geq -n_q$ for all $\gam \in \Lie \big(\Divf_{\Ct/C}^{0}\big)$. 
The bound $-n_q$ %for $\val_q \big( \Lie\big(\Divf_{\Ct/C}^{0}\big) \big)$ 
depends only on the singularity of $C$ at $q$. 
More precisely, $n_q$ satisfies 
$\fm_{\Ct,q}^{n_q+1} \subset \fm_{C,p}$ (cf.\ Lemma \ref{IDiv_Y/X^0 repr}). 
The number $n_q$ is related to the dimension of the affine part of $\Pic_C$, 
see \cite[Section 9.2, proof of Proposition 9]{BLR}. 
If $\sL$ is a sufficiently ample line bundle on $X$ 
(i.e.\ a sufficiently high power of an ample line bundle on $X$), 
one finds a family $T \subset \left|\sL\right|^{d-1}$ of Cartier curves $C$ 
whose proper transforms $C^Y$ intersect $\Es_{\eta}$ properly in points $p_C$ 
such that the set $\{p_C \,|\, C \in T\}$
contains an open dense subset $U_{\eta}$ of $\Es_{\eta}$. 
By upper semi-continuity of $\dim \Pic_C$ for the curves $C$ 
in $\left|\sL\right|^{d-1}$, we may assume 
that the sets $\val_{q_C} \big( \Lie\big(\Divf_{\Ct/C}^{0}\big) \big)$ 
for $q_C \in \mu^{-1}(p_C)$ admit a common bound $-n_{\eta}$ for all $C \in T$. 
Then for each $\del \in \Lie\big(\Divf_{Y/X}^{0}\big)$ 
with $\eta\in\Supp(\del)$ 
there is a curve $C \in T$ 
with $\val_{q_C}\big(\del\cut\Ct\big) = \val_{\eta}\left(\del\right)$, 
since this is an open dense property 
among the curves that intersect $\Es_{\eta}$. 
By definition, 
$\del \in \Lie\big(\Divf_{Y/X}^{0}\big)$ implies that 
$\del\cut\Ct \in \Lie\big(\Divf_{\Ct/C}^{0}\big)$. 
Then 
$\val_{\eta}\left(\del\right) = \val_{q_C}\big(\del\cut\Ct\big) \geq -n_{\eta}$. 
We obtain \; 
$\min \big(\val_{\eta}\big(\Lie\big(\Divf_{Y/X}^{0}\big)\big)\big) \geq -n_{\eta}$, 
i.e.\ the orders of poles of deformations in $\Lie\big(\Divf_{Y/X}^{0}\big)$ 
are bounded. %from below. 
Hence for all $\eta\in\SgY^{\htl}$ there exist $n_{\eta}$ such that 
\[  \Lie\left(\Divf_{Y/X}^{0}\right) \;\subset\; 
    \Gam\left(\left.\sO_Y\left(\sum_{\eta\in \SgY^{\htl}}n_{\eta}\,\Es_{\eta}\right)
              \right/\sO_Y\right) \;. 
\] 
As $Y$ is projective, the $k$-vector space on the right hand side 
is finite dimensional. 
\ePf 

\bRmk 
The construction of $\Divf_{Y/X}$ involves an intersection ranging over 
all Cartier curves in $X$, which makes this object hard to grasp. 
In fact, once a formal subgroup $\fmlE$ of $\Divf_Y^0$ 
containing $\Divf_{Y/X}^0$ is found, 
$\Divf_{Y/X}^0$ can be computed from one single curve: 
Let (for example) $\fmlE$ be the formal group defined by 
$\fmlE(k) = \{D \in \Divf_Y^0 \,|\, \Supp(D) \subset \SgY \}$, 
$\Lie\left(\fmlE\right) = 
 \Gam\big(\big.\sO_Y\big(\sum_{\eta\in \SgY^{\htl}}n_{\eta}\,\Es_{\eta}\big) 
 \big/\sO_Y\big)$, 
as in the proof of Proposition \ref{Div_Y/X^0}. 
Then there exists a Cartier curve $C$ in $X$ relative to $X_{\sing}$
such that 
\[  \Divf_{Y/X}^0 = 
    \left.\left( \_\cut\Ct \right)\right|_{\fmlE}^{-1} \Divf_{\Ct/C}^0 \;. 
\]
\eRmk 

\bPf[Indication of Proof] 
For each Cartier curve $C$ in $X$,  %consider 
define a subfunctor $\fmlF_C$ of $\fmlE$ by 
$ \fmlF_C := %\fmlE \isec 
  \big.\big( \_\cut\Ct \big)\big|_{\fmlE}^{-1} \Divf_{\Ct/C}^0  
$. 
Then $\fmlF_C$ is a formal group for each $C$, 
and it holds $\Divoyx = \bigcap_{C} \fmlF_C$. 
For each sequence $\{C_{\nu}\}$ of Cartier curves the formal groups 
$\fmlE_0 := \fmlE$, $\fmlE_{\nu+1} := \fmlE_{\nu} \isec \fmlF_{C_{\nu}}$ 
form a descending chain 
\[ \fmlE_0 \,\supset\, \fmlE_1 \,\supset\, \ldots \,\supset\, \fmlE_{\nu} 
   \,\supset\, \ldots \,\supset\, \bigcap_{C} \fmlF_{C} 
\] 
of formal subgroups of $\fmlE$. 
It is obvious that if $rD \in \Divf_{Y/X}^0(k)$ for some $D\in\Divf_Y(k)$ 
and $r\in\Zint\setminus\{0\}$, then $D \in \Divf_{Y/X}^0(k)$. 
Therefore $\Divoyx(k)$ is generated by a subset of a set of generators of 
$\fmlE(k)$, and this is a finitely generated free abelian group. 
Moreover, $\Lie\left(\fmlE\right)$ 
is a finite dimensional $k$-vector space. 
Hence the sequence $\{\fmlE_{\nu}\}$ becomes stationary. 
Thus there is a finite set of Cartier curves 
$C_1,\ldots,C_r$ such that 
$\Divf_{Y/X}^0 = \fmlF_{C_1} \isec \ldots \isec \fmlF_{C_r}$. 
Then each Cartier curve $C$ containing $C_1,\ldots,C_r$ 
gives the desired Cartier curve. 
\ePf

\subsection{The Category $\MrCH{X}$ } 
\label{subsec:Mr_CH_0(X)_deg0}

We keep the notation fixed at the beginning of Section 
\ref{sec:Rat-Maps CH_0(X)_deg0}: 
$X$ is a projective variety, 
$\pi:Y\lra X$ a projective resolution of singularities 
and $U\subset X_{\reg}$ a dense open subset. 

\bDef \label{MrCH(X)} 
$\MrCH{X}$ is a category of rational maps 
from $X$ to algebraic groups defined as follows: 
The objects of $\MrCH{X}$ are morphisms $\phe:U\lra G$ 
whose associated map on zero-cycles of degree zero 
  $\ZOo{U}  \lra  G(k)$, \; 
  $\sum l_i \, p_i  \lmt  \sum l_i \, \phe(p_i)$ 
factors through a homomorphism of groups $\CHOo{X}\lra G(k)$. 
\footnote{A category of rational maps to algebraic groups is defined 
already by its objects, according to Remark \ref{EquCatMr}. } \\ 
We refer to the objects of $\MrCH{x}$ as 
rational maps from $X$ to algebraic groups 
\emph{factoring through rational equivalence} or 
\emph{factoring through $\CHOo{X}$}. 
\eDef 

\bThm  \label{Equivalence} 
The category $\MrCH{X}$ of morphisms
from $U$ to algebraic groups factoring through $\CHOo{X}$
is equivalent to the category $\Mr_{\Divf_{Y/X}^{0}}$ of rational
maps from $Y$ to algebraic groups which induce a transformation of
formal groups to $\Divf_{Y/X}^{0}$ 
(see Proposition \ref{Div_Y/X^0} for the Definition of $\Divf_{Y/X}^{0}$). 
\eThm 

\bPf  
First notice that a rational map from $Y$ to an algebraic
group which induces a transformation to $\Divf_{Y/X}^{0}$ is necessarily
regular on $U$, since all $D\inn\Divf_{Y/X}^{0}$ have support only
on $Y\setminus U$. Then according to Definition \ref{CH_0(X)_deg0}
and Definition \ref{Mr_sW} the task is to show that for a morphism
$\phe:U\lra G$ from $U$ to an algebraic group $G$
with canonical decomposition 
$0\ra L\ra G\ra A\ra0$
the following conditions are equivalent: 

\begin{tabular}{rl}
(i)&
$\phe\left(\dv\left(f\right)_{C}\right)=0\qquad\qquad
\forall\left(C,f\right)\in\fR_{0}\left(X,U\right)$ \tabularnewline
(ii)&
$\dv_{\splG}\left(\phe_{Y,\lambda}\right)\inn\Divf_{Y/X}^{0}\qquad\qquad
\forall\lambda\inn L^{\vee}$\tabularnewline
\end{tabular} \\
where $\phe_{Y,\lambda}$ is the induced section of the $\splG$-bundle
$\lambda_{*}G_{Y}$ over $Y$ introduced in Subsection 
\ref{sub:Categories-of-Rational}.
The principal $L$-bundle $G$ is a direct sum of $\splG$-bundles $\lambda_{*}G$
over $A$, $\lambda\inn\Ld$; let $\phe_{\lambda}:U\lra\lambda_{*}G$
be the induced morphisms. Then condition (i) is equivalent to 

\begin{tabular}{rl}
(i')&
$\phe_{\lambda}\left(\dv\left(f\right)_{C}\right)=0\qquad\qquad
\forall\lambda\inn L^{\vee},\;\forall\left(C,f\right)\in\fR_{0}\left(X,U\right)$. 
\tabularnewline
\end{tabular} \\
Hence it comes down to show that for all $\lambda\inn\Ld$ the following
conditions are equivalent: 

\begin{tabular}{rl}
(j)&
$\phe_{\lambda}\left(\dv\left(f\right)_{C}\right)=0\qquad\qquad
\forall\left(C,f\right)\in\fR_{0}\left(X,U\right)$, \tabularnewline
(jj)&
$\dv_{\splG}\left(\phe_{Y,\lambda}\right)\inn\Divf_{Y/X}^{0}$\;.\tabularnewline
\end{tabular} \\
This is the content of Lemma \ref{L-bdl_factors_CH_0} below. 
\ePf 

\bLem \label{L-bdl_factors_CH_0} 
Let $\phe_{\lambda}:U\lra G_{\lambda}$
be a morphism from $U$ to an algebraic group 
$G_{\lambda}\in\Ext\left(A,\splG\right)$,  
i.e.\ $G_{\lambda}$ is a $\splG$-bundle over an abelian variety $A$. 
Then the following conditions are equivalent: 

\begin{tabular}{rl}
(i)&
$\phe_{\lambda}\left(\dv\left(f\right)_{C}\right)=0\qquad\qquad
\forall\left(C,f\right)\in\fR_{0}\left(X,U\right)$, \tabularnewline
(ii)&
$\dv_{\splG}\left(\phe_{Y,\lambda}\right)\inn\Divf_{Y/X}^{0}$ . \tabularnewline
\end{tabular} \eLem 

\bPf  
Let $C$ be a Cartier curve in
$X$ relative to $X\setminus U$, and let $\nu:\Ct\lra C$
be its normalization. 
In the case $\splG=\Gm$ Lemma \ref{Gm-bdl_factors_CH_0} and in the case 
$\splG=\Ga$ Lemma \ref{Ga-bdl_factors_CH_0} assert that the following
conditions are equivalent:

\begin{tabular}{rl}
(j)&
$\phe_{\lambda}|_{C}\left(\dv\left(f\right)\right)=0\qquad\qquad
\forall f\in\K\left(C,C\cap U\right)^{*}$,\tabularnewline
(jj)&
$\nu_{*}\left(\dv_{\splG}\left(\phe_{\lambda}|_{C}\right)\right)=0$\;. \tabularnewline
\end{tabular} \\
We have $\dv_{\splG}\left(\phe_{\lambda}|_{C}\right) = 
\dv_{\splG}\left(\phe_{Y,\lambda}\right) \cut \Ct$, where
$\_\cut \Ct:\Decf_{Y,\Ct}\lra\Divf_{\Ct}$ is the pull-back of Cartier 
divisors from $Y$ to $\Ct$ 
($\see$Definition \ref{Dec_Y,V}, Proposition \ref{_.V}). 
Using the equivalence (j)$\Llra$(jj) above, condition (i) is equivalent to 

\begin{tabular}{rl}
(i')&
$\left(\nu_C\right)_{*}
\left(\dv_{\splG}\left(\phe_{Y,\lambda}\right) \cut \Ct \right) = 0
\qquad \forall$\,Cartier curves $C$ relative to $X\setminus U$. \\
\end{tabular} \\
Conditions (i') and (ii) are equivalent by definition of $\Divf_{Y/X}^{0}$
($\see$Proposition \ref{Div_Y/X^0}), taking into account that
$\dv_{\splG}\left(\phe_{Y,\lambda}\right)\inn\Divf_{Y}^{0}$ by Proposition
\ref{induced Trafo}. 
\ePf 

\bLem \label{Gm-bdl_factors_CH_0} Let $C$ be a projective curve
and $\nu:\Ct\lra C$ its normalization. Let 
$\psi:C\dra G_{\mu}$ be a rational map from $C$ 
to an algebraic group $G_{\mu}\in\Ext\left(A,\Gm\right)$, 
i.e.\ $G_{\mu}$ is a $\Gm$-bundle over an abelian variety $A$.
Suppose that $\psi$ is regular on a dense open subset $U_{C}\subset C_{\reg}$,
which we identify with its preimage in $\Ct$. Then the following
conditions are equivalent: 

\begin{tabular}{rl}
(i)&
$\psi\left(\dv\left(f\right)\right)=0\qquad\qquad\qquad\quad\; 
\forall f\in\K\left(C,U_{C}\right)^{*}$\;, \tabularnewline
(ii)&
$\left(f\circ\nu\right)\left(\dv_{\Gm}\left(\psi\right)\right)=0\qquad\qquad
\forall f\in\K\left(C,U_{C}\right)^{*}$\;, \tabularnewline
(iii)&
$\nu_{*}\left(\dv_{\Gm}\left(\psi\right)\right)=0$\;. \tabularnewline
\end{tabular} \eLem 

\bPf  
(i)$\Llra$(ii) We show that for all 
$f\in\K\left(C,U_{C}\right)^{*}$ it holds 
\[ \psi\left(\dv\left(f\right)\right)=
   \left(f\circ\nu\right)\left(\dv_{\Gm}\left(\psi\right)\right) \;. 
\] 
Let $f\in\K\left(C,U_{C}\right)^{*}$. 
Write $\widetilde{f} := \nu^{\fis}f = f\circ\nu$. 
Set $S:=\Ct\setminus U_{C}$. For each $s\in S$ let 
$\Phi_{s}:U_{s}\times\Gm\lra G_{\mu}$ 
be a local trivialization of the induced $\Gm$-bundle over $\Ct$
in a neighbourhood $U_{s}\ni s$. 
Notice that 
$\val_{p}\left(\psi\right):=\val_{p}\left(\left[\psi\right]_{\Phi_{p}}\right)$
is independent of the local trivialization. 
Since $f\in\K\left(C,U_{C}\right)^{*}$, 
we have $f\in\sO_{C,s}^{*}$ for all $s\in S$ and hence 
$\dv\left(f\right)\cap S=\varnothing$.
Then by Lemma \ref{normlzd} it holds 
\[ \psi\left(\dv\left(f\right)\right) 
   = \left(\psi\circ\nu\right)\left(\dv\left(\widetilde{f}\right)\right) 
   = \prod_{c\notin S}\psi(c)^{\val_{c}\left(\widetilde{f}\right)} \;. 
\] 
The defining properties of a local symbol from Definition \ref{locSymbol} 
imply 
\[ \prod_{c\notin S}\psi(c)^{\val_{c}\left(\widetilde{f}\right)} 
   = \prod_{c\notin S}\left(\psi,\widetilde{f}\right)_{c} 
   = \prod_{s\in S}\left(\psi,\widetilde{f}\right)_{s}^{-1} \;. 
\] 
According to Lemma \ref{locSymb_triv} and the explicit description
from Example \ref{locSymb_Gm} of local symbols for rational maps
to $\Gm$ this is equal to 
\[ \prod_{s\in S}\left(\psi,\widetilde{f}\right)_{s}^{-1} 
   = \prod_{s\in S}\left(\left[\psi\right]_{\Phi_{s}},\widetilde{f}\right)_{s}^{-1} 
   = \prod_{s\in S}\left(\widetilde{f},\left[\psi\right]_{\Phi_{s}}\right)_{s} \;. 
\] 
Finally, using again the defining properties of a local symbol 
from Definition \ref{locSymbol}, we obtain 
\[ \prod_{s\in S}\left(\widetilde{f},\left[\psi\right]_{\Phi_{s}}\right)_{s} 
   = \prod_{p\notin\Supp\left(\dv\left(\widetilde{f}\right)\right)} 
     \widetilde{f}(p)^{\val_{p}\left(\psi\right)} 
   = \widetilde{f}\left(\dv_{\Gm}\left(\psi\right)\right) \;. 
\] 

(ii)$\Llra$(iii) The implication (iii)$\Lra$(ii) 
is clear. For the converse direction first observe that the support of 
$\dv_{\Gm}(\psi)$ lies necessarily in $\Ct\setminus U_C$, since $\psi$ is 
regular on $U_C$. For each $s \in C\setminus U_C$ there is a 
rational function $f_s \in \K\left(C,U_{C}\right)^{*}$ such that 
$f(s)=t\in\Gm\setminus\{1\}$ and $f(z)=1$ for all 
$z \in C\setminus (U_C\cup\{s\})$ by the approximation theorem. 
Then $\left(f_s\circ\nu\right)\left(\dv_{\Gm}\left(\psi\right)\right)=0$ 
if and only if $\nu_{*}\left(\dv_{\Gm}\left(\psi\right)|_{\nu^{-1}(s)}\right)=0$, 
where $\dv_{\Gm}\left(\psi\right)|_{\nu^{-1}(s)}$ is the part of 
$\dv_{\Gm}\left(\psi\right)$ which has support on $\nu^{-1}(s)$. 
As this is true for all $s \in C\setminus U_C$, it shows the implication 
(ii)$\Lra$(iii). 
\ePf 

\bLem \label{Ga-bdl_factors_CH_0} Let $C$ be a projective curve
and $\nu:\Ct\lra C$ its normalization. 
Let $\psi:C\dra G_{\alpha}$ be rational map from $C$ 
to an algebraic group $G_{\alpha}\in\Ext\left(A,\Ga\right)$, 
i.e.\ $G_{\alpha}$ is a $\Ga$-bundle over an abelian variety $A$.
Suppose that $\psi$ is regular on a dense open subset $U_{C}\subset C_{\reg}$,
which we identify with its preimage in $\Ct$. Then the following
conditions are equivalent: 

\begin{tabular}{rl}
(i)&
$\psi\left(\dv\left(f\right)\right)=0\qquad\qquad\qquad\quad
 \forall f\in\K\left(C,U_{C}\right)^{*}$\;, \tabularnewline
(ii)&
$\sum_{q\ra p}\Res_{q}\left(\psi\;\der g\right)=0
 \qquad\qquad\forall g\in\Oc_{C,p},\;\forall p\in C$\;, \tabularnewline
(iii)&
$\nu_{*}\left(\dv_{\Ga}\left(\psi\right)\right)=0$\;. \tabularnewline
\end{tabular} \eLem 

\bPf  (i)$\Llra$(ii) Let $f\in\K\left(C,U_{C}\right)^{*}$.
We will identify $f$ with $f':=\nu^{\fis}f$. %, if no confusion is likely. 
Set $S:=\Ct\setminus U_{C}$. %and $S:=C\setminus U_{C}$. 
For each $s\in S$ let $\Phi_{s}:U_{s}\times\Ga\lra G_{\alpha}$
be a local trivialization of the induced $\Ga$-bundle over $\Ct$
in a neighbourhood $U_{s}\ni s$. 
Notice that for each $\omega\in\Omega_{\Ct}$ 
which is regular at $q\in\Ct$ the expression 
$\Res_{q}\left(\psi\;\omega\right):=
\Res_{q}\left(\left[\psi\right]_{\Phi_{q}}\omega\right)$
is independent of the local trivialization. 
Then by Lemma \ref{normlzd} it holds 
\[ \psi\left(\dv\left(f\right)\right) 
   = \left(\psi\circ\nu\right)\left(\dv\left(f'\right)\right) 
   = \sum_{c\notin S}\val_{c}\left(f'\right)\;\psi(c) \;. 
\] 
The defining properties of a local symbol from Definition \ref{locSymbol} 
imply 
\[ \sum_{c\notin S}\val_{c}\left(f'\right)\;\psi(c) 
   = \sum_{c\notin S}\left(\psi,f'\right)_{c} 
   = -\sum_{s\in S}\left(\psi,f'\right)_{s} \;. 
\] 
According to Lemma \ref{locSymb_triv} and 
the explicit description from Example \ref{locSymb_Ga} of local symbols
for rational maps to $\Ga$ we obtain
\[ -\sum_{s\in S}\left(\psi,f'\right)_{s} 
   = -\sum_{s\in S}\left(\left[\psi\right]_{\Phi_{s}},f'\right)_{s} 
   = -\sum_{s\in S}\Res_{s}\left(\psi\;\der f'/f'\right) \;. 
\] 
Now $\der f/f=\der\log f$. 
More precisely, if $f=\alp(1+h)\in\sO_{C,p}^*$ with $\alp\in k^*$ 
and $h \in \fm_{C,p}$, then 
$\der f/f=\der (1+h)/(1+h)=\der\log(1+h)$ 
and $\log:1+\mc_{C,p}\overset{\sim}\lra\mc_{C,p}$ is well-defined. 
Thus the implication (ii)$\Lra$(i) is clear. 

For the converse, we first show that 
for each $p\in\nu(S)$, each $g\in\Oc_{C,p}$ 
and each effective divisor $\mdll$ %=\sum_{s\in S} m_s\, s$ 
supported on $S$ 
there is a rational function $f_{p}\in\K\left(C,U_C\right)^{*}$
such that $\der\log f_{p} \equiv \der g\mod\mdll$ at $\nu^{-1}(p)$ and 
$f_{p}\equiv1\mod\mdll$ at $s\in S\setminus\nu^{-1}(p)$. 
Since $\im \left(\Oc_{C,p}\overset{\der}\lra\Omega_{\Oc_{C,p}}\right)= 
\im \left(\mc_{C,p}\overset{\der}\lra\Omega_{\Oc_{C,p}}\right)$, 
we may assume $g\in\mc_{C,p}$. 
According to the approximation theorem, 
for $\mdl=\sum_{s\in S} n_s\, s$ there is $f_p\in\sK_{\Ct}$ 
such that $f_p \equiv \exp g \mod\mdl$ at $\nu^{-1}(p)$ and 
$f_{p}\equiv1\mod\mdl$ at $s\in S\setminus\nu^{-1}(p)$. 
One sees that $f_{p}\in\K\left(C,U_C\right)^{*}$. 
In particular, there is $h\in\mc_{C,p}$ 
with $h\in\mc_{\Ct,q}^{n_q}$ for each $q\ra p$ 
such that \;$\exp g = f_p+h = f_p(1+f_p^{-1}h)$. 
Then $g = \log f_p + \log(1+f_p^{-1}h)$, 
where $\log(1+f_p^{-1}h) \in \mc_{\Ct,q}^{n_q}$, 
since $f_p^{-1}h \in \mc_{\Ct,q}^{n_q}$ 
and $\log(1+\mc_{\Ct,q}^{n}) = \mc_{\Ct,q}^{n}$ for $n \geq 1$. 
This yields $\der\log f_{p} \equiv \der g\mod\mdll$ at $\nu^{-1}(p)$ 
if $\mdll=\sum_{s\in S} m_s\, s$ and $\mdl=\sum_{s\in S} (m_s+1)\, s$. 

Choosing $\mdll=\sum_{s\in S} m_s\, s$ large enough, 
i.e.\ $m_s$ larger than the pole order of $\psi$ at $s$, yields that 
$\Res_s\left(\psi\;\der f_p/f_p\right)=0$ for all $s\in S\setminus\nu^{-1}(p)$, 
as $\der f_p/f_p$ has a zero of order $\geq m_s -1$ at 
$s\in S\setminus\nu^{-1}(p)$. 
Hence $\psi\left(\dv\left(f_{p}\right)\right) = 0$ if and only if 
$\sum_{q\ra p}\Res_{p}\left(\psi\;\der f_{p}/f_{p}\right)=
 \sum_{q\ra p}\Res_{p}\left(\psi\;\der g\right)=0$. 
%Thus $\psi\left(\dv\left(f\right)_{C}\right)=0$
%for all $f\in\K\left(C,U_{C}\right)^{*}$ if and only if 
%$\sum_{q\ra p}\Res_{p}\left(\psi\;\der g\right)=0$
%for all $g\in\Oc_{C,p}$, $p\in\nu(S)$. 
It remains to remark that $\Res_{c}\left(\psi\;\der h\right)=0$ for all 
$h\in\Oc_{\Ct,c}\supset\Oc_{C,\nu\left(c\right)}$,
$c\in U_{C}$, since $\psi$ and $\der h$ are both regular at $c$. 

(ii)$\Llra$(iii) Let $q\in\Ct$. Then 
$\sum_{q\ra p}\Res_{p}\left(\psi\;\der g\right)=0$
for all $g\in\Oc_{C,p}$ is equivalent to the condition that the image 
$\sum_{q\ra p}\left[\dv_{\Ga}\left(\psi\right)\right]_{q}$ of 
$\dv_{\Ga}\left(\psi\right)$ in 
$\dsum_{q\ra p}\Hom_{k(q)}^{\cont}\left(\mc_{\Ct,q},k(q)\right)$ 
vanishes on $\mc_{C,p}$, 
by construction (%$\see$proof of 
Proposition \ref{fml}), 
which says 
$0=\sum_{q\ra p}\left[\dv_{\Ga}\left(\psi\right)\right]_{q}\circ\widehat{\nu^{\fis}}
\in\dsum_{q\ra p}\Hom_{k(q)}^{\cont}\left(\mc_{C,\nu(q)},k(q)\right)$. 
This is true for all $p\in C$ if and only if 
$\nu_{*}\left(\dv_{\Ga}\left(\psi\right)\right)=0$
by definition of the push-forward for formal infinitesimal divisors
($\see$Proposition \ref{phi_*}). 
\ePf 

\bLem \label{normlzd} 
Let $C$ be a Cartier curve in $X$ relative to $X\setminus U$ 
and $\nu:\Ct\lra C$ its normalization. 
If $\psi:C\cap U\lra G$ is a morphism
from $C\cap U$ to an algebraic group $G$, then for each 
$f\in\K\left(C,C\cap U\right)^*$ it holds \[
\psi\left(\dv\left(f\right)_{C}\right)=
\left(\psi\circ\nu\right)\left(\dv\left(\nu^{\fis}f\right)_{\Ct}\right) \;. 
\] 
\eLem 

\bPf Follows immediately from Definition \ref{div(f)_C}. 
\ePf

\subsection{Universal Regular Quotient} 
\label{subsec:Universal_Regular_Quotient} 

The results obtained up to now provide the necessary foundations for a 
description of the universal regular quotient and its dual, which was the 
initial intention of this work. 

\subsubsection*{Existence and Construction} 
\label{subsubsec:Exist_UnivRegQuot} 

The universal regular quotient $\Alb\left(X\right)$ of a (singular)
projective variety $X$ is by definition (see \cite{ESV}) 
the universal object for the category $\MrCH{X}$ 
of morphisms from $U \subset X_\reg$ 
factoring through $\CHOo{X}$ ($\see$Definition \ref{MrCH(X)}). 
%which factor through a homomorphism of groups $\CHOo{X}\lra G$. 
In Theorem \ref{Equivalence} we have seen that this category is equivalent 
to the category $\Mr_{\Divf_{Y/X}^0}$ of rational maps from 
a projective resolution of singularities $Y$ for $X$ to algebraic groups 
which induce a transformation to the formal group $\Divf_{Y/X}^0$. 
Now Theorem \ref{Exist univObj} implies the existence 
of a universal object $\Alb_{\Divf_{Y/X}^0}(Y)$ for this category, which was 
constructed ($\see$Remark \ref{Alb_constr}) as the dual 1-motive of 
$\left[ \Divf_{Y/X}^0 \lra \Picf_Y^0 \right]$. 
As $\Alb\left(X\right) = \Alb_{\Divf_{Y/X}^0}(Y)$, this gives the 
existence and an explicit construction of the universal regular quotient, 
as well as a description of its dual. The proof of Theorem \ref{main_result} 
is thus complete. 

\subsubsection*{Functoriality} 
\label{subsubsec:Fctr_UnivRegQuot} 

Let $X$, $V$ be projective varieties whose normalizations 
$\Xt$, $\Vt$ are regular. 
We analyze  whether a morphism $\sigma: V \lra X$ induces 
a homomorphism of algebraic groups $\Alb \lp V \rp \lra \Alb \lp X \rp$. 

As the functoriality of the universal objects $\Alb_{\fmlG}(Y)$, where $Y$ 
is regular and $\fmlG \subset \Divf_Y^0$ is a formal group, has already been 
treated in Proposition \ref{alb_F(sigma)}, %on page \pageref{alb(sigma)}, 
we will reduce the problem to this case. 
Therefore it obliges to show under which assumptions the 
following conditions hold: 

\begin{tabular}{rl} 
($\alp$)  & A morphism $\sigma: V \lra X$ induces 
           a morphism $\sigt: \Vt \lra \Xt$. \\
($\beta$) & The pull-back of relative Cartier divisors maps 
           $\Divf_{\Xt/X}^0$ to $\Divf_{\Vt/V}^0$. \\
\end{tabular} \\ 
For this purpose we introduce the following notion, analogue to Definition 
\ref{Cartier-curve} (keeping the notation fixed at the beginning of 
this Section \ref{sec:Rat-Maps CH_0(X)_deg0}): 

\bDef \label{Cartier-morphism} 
A \emph{Cartier subvariety in $X$ relative to $X\setminus U$} 
is a subvariety $V \subset X$ satisfying \\ 
\begin{tabular}{rl}
(a) & $V$ is equi-dimensional. 
      \tabularnewline
(b) & No component of $V$ is contained in $X \setminus U$. 
      \tabularnewline
(c) & If $p \in V\setminus U$, 
      the ideal of $V$ in $\sO_{X,p}$ 
      is generated by a regular sequence.  \tabularnewline
\end{tabular} 
\eDef 

\bRmk 
%The terminology \emph{Cartier subvariety} might be confusing, because 
A Cartier subvariety $V$ in $X$ relative to $X\setminus U$ 
\emph{in codimension one} needs not to be a Cartier divisor 
on the whole of $X$. 
Point (c) of Definition \ref{Cartier-morphism} implies that $V$ 
is a locally principal divisor in a neighbourhood of $X\setminus U$. 
\eRmk 

\bPrp \label{Pic-1-mot_functorial} 
Let $V \subset X$ be a Cartier subvariety relative to $X\setminus U$. 
Then the pull-back of relative Cartier divisors and of line bundles 
induces a transformation of 1-motives 
\[ \left[ \begin{array}{c} 
          \Divf_{\Vt/V}^0 \\ \downarrow \\ \Picf_{\Vt}^0 \\
          \end{array} 
   \right]
   \lla 
   \left[ \begin{array}{c} 
          \Divf_{\Xt/X}^0 \\ \downarrow \\ \Picf_{\Xt}^0 \\
          \end{array} 
   \right] \;. 
\] 
\ePrp 

\bPf It suffices to verify the conditions ($\alp$) and ($\beta$) 
mentioned above. \\ 
%at the beginning of this Subsubsection \ref{subsubsec:Fctr_UnivRegQuot}. 
As no irreducible component of $V$ is contained in $X\setminus U$ 
by condition (b) of Definition \ref{Cartier-morphism}, 
the base change $V \tms_X \Xt =: V_{\Xt} \lra V$ of 
$\Xt \lra X$ is biregular for each irreducible component of $V_{\Xt}$, 
and there is exactly one irreducible component of $V_{\Xt}$ 
lying over each irreducible component of $V$. 
%a morphism of degree 1 ($\see$Definition \ref{deg1}) 
Thus the normalization $\Vt \lra V$ factors through $V_{\Xt} \lra V$. 
We obtain a commutative diagram 
\[ \xymatrix{ \Vt \ar@/^1pc/[drr] \ar@/_1pc/[ddr] \ar[dr] & & \\
              & V_{\Xt} \ar[r] \ar[d] & \Xt \ar[d] \\ 
              & V \ar[r] & X \\ 
            } 
\] 
The morphism $\Vt \lra \Xt$ induces a pull-back of families of line bundles 
$\Picf_{\Xt}^0 \lra \Picf_{\Vt}^0$ and a pull-back of relative Cartier divisors 
$\Divf_{\Xt/X}^0 \lra \Divf_{\Vt}^0$, since no component of $V$ is contained 
in $\Supp \lp \Divf_{\Xt/X}^0 \rp$. 
A Cartier curve in $V$ relative to $U_V = V \tms_X U$ is also 
a Cartier curve in $X$ relative to $X\setminus U$. Therefore 
the definition of $\Divf_{\Vt/V}^0$ ($\see$Proposition \ref{Div_Y/X^0})
implies that the image of $\Divf_{\Xt/X}^0$ under pull-back $\_\cut \Vt$ 
lies actually in $\Divf_{\Vt/V}^0$. 
This gives a commutative diagram of 
natural transformations of functors 
\[ \xymatrix{ \Divf_{\Vt/V}^0 \ar[d] & \Divf_{\Xt/X}^0 \ar[d] \ar[l] \\ 
              \Picf_{\Vt}^0 & \Picf_{\Xt}^0 \ar[l] \\ 
            } %\;. 
\] 
\ePf 
\vspace{\vs} 

Dualization of 1-motives yields the following functoriality 
of the universal regular quotient: 

\bPrp \label{UnivRegQuot_functorial} 
Let $\iota:V \subset X$ be a Cartier subvariety relative to $X\setminus U$. 
Then $\iota$ induces a homomorphism of algebraic groups 
\[ \Alb(\iota): \Alb \lp V \rp \lra \Alb \lp X \rp  \;. 
\] 
\ePrp

\vspace{\vs}

{\scshape
\begin{flushright}
\begin{tabular}{l}
Universit\"at Duisburg-Essen \\
FB6 Mathematik, Campus Essen \\
45117 Essen \\
Germany \\
{\upshape e-mail: henrik.russell@uni-due.de}\\
\end{tabular}
\end{flushright}
}

\end{document}